# A MICROSCOPIC PROBABILISTIC DESCRIPTION OF A LOCALLY REGULATED POPULATION AND MACROSCOPIC APPROXIMATIONS

By Nicolas Fournier and Sylvie Méléard

*Institut Elie Cartan and Université Paris 10*

We consider a discrete model that describes a locally regulated spatial population with mortality selection. This model was studied in parallel by Bolker and Pacala and Dieckmann, Law and Murrell. We first generalize this model by adding spatial dependence. Then we give a pathwise description in terms of Poisson point measures. We show that different normalizations may lead to different macroscopic approximations of this model. The first approximation is deterministic and gives a rigorous sense to the *number density*. The second approximation is a superprocess previously studied by Etheridge. Finally, we study in specific cases the long time behavior of the system and of its deterministic approximation.

**1. Introduction.** We consider a spatial ecological system that consists of motionless individuals (such as *plants*). Individuals are characterized by their location. We assume that each plant produces seeds at a given rate. When a seed is born, it immediately disperses from its *mother* and becomes a mature plant. We also assume that plants are subjected to *mortality selection*. That is, each plant dies at a rate that depends on the local population density. All these events occur randomly in continuous time. This model was introduced by Bolker and Pacala [2] and Dieckmann and Law [9]. To study the system, Bolker and Pacala derived approximations for the time evolution of the moments (mean and spatial covariance) of the population distribution. In the present article, we wish to give a rigorous definition of the underlying *microscopic* stochastic process and rewrite rigorously the moment equations of [2], then to derive some tractable macroscopic approximations, and finally to study the long time behavior of the stochastic process and its approxi-mations. Unfortunately, we obtained only partial results concerning the last point.

In Section 2, we describe the Bolker–Pacala–Dieckmann–Law (BPDL) process in detail. In fact, we generalize the model slightly by adding a spatial dependence in all the rates. Then we give a pathwise representation of the system in terms of Poisson point measures. We also produce a numerical algorithm to simulate the BPDL process. Section 3 is devoted to existence and uniqueness. We also show some martingale properties of the BPDL process. In Section 4, we find the mean equation that Bolker and Pacala [2] intuitively obtained. We also give a rigorous sense to the covariance terms formally defined in [2] or [9], [4] and [10]. Section 5 is concerned with macroscopic approximations of the BPDL process. We first show that, conveniently normalized, the BPDL process converges to the solution of a deterministic nonlinear integrodifferential equation. We propose this as a rigorous interpretation of the *density number*, often introduced by biologists without a proper definition. We also show that with another normalization, the BPDL process converges to the superprocess version of the BPDL model introduced and studied by Etheridge [6]. We give partial results about extinction and survival for the BPDL process in Section 6. In Section 7, we study the convergence to equilibrium of the deterministic approximation. We obtain only some partial results. We next show that in the *detailed balance case* to be specified later on, there exists a nontrival steady state for the BPDL process. We conclude the article with some simulations.

**2. The model.** Let us first describe the model in detail.

2.1. *Definition of the parameters and heuristics.* The plants are supposed to be motionless and characterized by their spatial location. We assume that the spatial domain is the closure $\bar{\mathcal{X}}$ of an open connected subset $\mathcal{X}$ of $\mathbb{R}^d$, for some $d \geq 1$. We denote by $M_F(\bar{\mathcal{X}})$ [resp. $\mathcal{P}(\bar{\mathcal{X}})$] the set of finite nonnegative measures (resp. probability measures) on $\bar{\mathcal{X}}$. Let also $\mathcal{M}$ be the subset of $M_F(\bar{\mathcal{X}})$ that consists of all finite point measures:

$$(2.1) \qquad \mathcal{M} = \left\{ \sum_{i=1}^n \delta_{x_i}, n \geq 0, x_1, \ldots, x_n \in \bar{\mathcal{X}} \right\}.$$

Here and below, $\delta_x$ denotes the Dirac mass at $x$. For any $m = \sum_{i=1}^n \delta_{x_i} \in \mathcal{M}$, any measurable function $f$ on $\bar{\mathcal{X}}$, we set $\langle m, f \rangle = \int_{\bar{\mathcal{X}}} f \, dm = \sum_{i=1}^n f(x_i)$.

NOTATION 2.1. For all $x$ in $\bar{\mathcal{X}}$, we introduce the following quantities:

(i) $\mu(x) \in [0, \infty)$ is the rate of "intrinsic" death of plants located at $x$,
(ii) $\gamma(x) \in [0, \infty)$ is the rate of seed production of plants located at $x$,

(iii) $D(x, dz)$ is the dispersion law of the seeds of plants located at $x$. It is assumed to satisfy, for each $x \in \bar{\mathcal{X}}$,

$$\int_{z \in \mathbb{R}^d, x+z \in \bar{\mathcal{X}}} D(x, dz) = 1 \quad \text{and} \quad \int_{z \in \mathbb{R}^d, x+z \notin \bar{\mathcal{X}}} D(x, dz) = 0.$$

(iv) $\alpha(x) \in [0, \infty)$ is the rate of interaction of plants located at $x$.

(v) For $x$, $y$ in $\bar{\mathcal{X}}$, $U(x,y) = U(y,x) \in [0, \infty)$ is the competition kernel.

The competition kernel $U(x,y)$ describes the strength of competition between plants located at $x$ and $y$.

We aim to study the stochastic process $\nu_t$, taking its values in $\mathcal{M}$ and describing the *distribution* of plants at time $t$. We write

(2.2) $$\nu_t = \sum_{i=1}^{I(t)} \delta_{X_t^i},$$

where $I(t) \in \mathbb{N}$ stands for the number of plants alive at time $t$ and $X_t^1, \ldots, X_t^{I(t)}$ describe their locations (in $\bar{\mathcal{X}}$). The supposed dynamics for this population can be roughly summarized as follows:

(i) At time $t = 0$, we have a (possibly random) distribution $\nu_0 \in \mathcal{M}$.

(ii) Each plant (located at some $x \in \bar{\mathcal{X}}$) has three independent exponential clocks: a *seed production* clock with parameter $\gamma(x)$, a *natural death* clock with parameter $\mu(x)$ and a *competition mortality* clock with parameter $\alpha(x) \sum_{i=1}^{I(t)} U(x, X_t^i)$.

(iii) If one of the two *death* clocks of a plant rings, then this plant disappears.

(iv) If the *seed production* clock of a plant (located at some $x \in \bar{\mathcal{X}}$) rings, then it produces a seed. This seed immediately becomes a mature plant. Its location is given by $y = x + z$, where $z$ is randomly chosen according to the dispersion law $D(x, dz)$.

In [2], $\gamma$, $\mu$, $\alpha$ and $D$ were assumed to be space-independent. Our generalization might allow us to take into account external effects such as landscape, resource distribution and so forth. Note also that assuming that all these clocks are exponentially distributed allows us to reset all the clocks to 0 each time one clock rings.

We wish to describe the system by the evolution in time of the empirical measure $\nu_t$. More precisely, we are looking for an $\mathcal{M}$-valued Markov process $(\nu_t)_{t \geq 0}$ with infinitesimal generator $L$, defined for a large class of functions $\phi$ from $\mathcal{M}$ into $\mathbb{R}$, for all $\nu \in \mathcal{M}$, by

(2.3)
$$L\phi(\nu) = \int_{\bar{\mathcal{X}}} \nu(dx) \int_{\mathbb{R}^d} D(x, dz)[\phi(\nu + \delta_{x+z}) - \phi(\nu)]\gamma(x)$$
$$+ \int_{\bar{\mathcal{X}}} \nu(dx)[\phi(\nu - \delta_x) - \phi(\nu)]\left\{\mu(x) + \alpha(x) \int_{\bar{\mathcal{X}}} \nu(dy) U(x,y)\right\}.$$

The first term is linear (in $\nu$) and describes the seed production and dispersal phenomenon. The second term is nonlinear and describes death due to age or competition. This infinitesimal generator can be compared with formula (3) in [2], page 182.

2.2. *Description in terms of Poisson measures.* We now give a pathwise description of the $\mathcal{M}$-valued stochastic process $(\nu_t)_{t\geq 0}$. To this end, we use Poisson point measures. For the sake of simplicity, we assume that the spatial dependence of all the parameters is bounded in some sense.

ASSUMPTION A. There exist some constants $\bar{\alpha}$, $\bar{\gamma}$ and $\bar{\mu}$ such that, for all $x \in \bar{\mathcal{X}}$,

$$(2.4) \qquad \alpha(x) \leq \bar{\alpha}, \qquad \gamma(x) \leq \bar{\gamma}, \qquad \mu(x) \leq \bar{\mu}.$$

There exist a constant $C > 0$ and a probability density $\tilde{D}$ on $\mathbb{R}^d$ such that, for all $x \in \bar{\mathcal{X}}$,

$$(2.5) \qquad D(x,dz) = D(x,z)\,dz \qquad \text{with } D(x,z) \leq C\tilde{D}(z).$$

The competition kernel $U$ is bounded by some constant $\bar{U}$.

We also introduce the following notation.

NOTATION 2.2. Let $\mathbb{N}^* = \mathbb{N} \setminus \{0\}$. Let $H = (H^1, \ldots, H^k, \ldots) : \mathcal{M} \mapsto (\mathbb{R}^d)^{\mathbb{N}^*}$ be defined by

$$(2.6) \qquad H\left(\sum_{i=1}^n \delta_{x_i}\right) = (x_{\sigma(1)}, \ldots, x_{\sigma(n)}, 0, \ldots, 0, \ldots),$$

where $x_{\sigma(1)} \preccurlyeq \cdots \preccurlyeq x_{\sigma(n)}$ for some arbitrary order $\preccurlyeq$ on $\mathbb{R}^d$ (one may, e.g., choose the lexicographic order).

This function $H$ allows us to overcome the following (purely notational) problem: Assume that a population of plants is described by a point measure $\nu \in \mathcal{M}$. Choosing a plant uniformly among all plants consists of choosing $i$ uniformly in $\{1, \ldots, \langle \nu, 1 \rangle\}$, and then choosing the plant *number i* (from the arbitrary order point of view). The location of such a plant is thus $H^i(\nu)$.

NOTATION 2.3. We consider the path space $\mathcal{T} \subset \mathbb{D}([0,\infty), M_F(\bar{\mathcal{X}}))$ defined by

$$(2.7) \quad \mathcal{T} = \left\{ (\nu_t)_{t\geq 0} \,\middle/\, \begin{array}{l} \forall t \geq 0, \nu_t \in \mathcal{M}, \text{ and } \exists\, 0 = t_0 < t_1 < t_2 < \cdots, \\ \lim_n t_n = \infty \text{ and } \nu_t = \nu_{t_i}\ \forall t \in [t_i, t_{i+1}) \end{array} \right\}.$$

Note that for $(\nu_t)_{t\geq 0} \in \mathcal{T}$, and $t > 0$ we can define $\nu_{t-}$ in the following way: If $t \notin \bigcup_i \{t_i\}$, $\nu_{t-} = \nu_t$, while if $t = t_i$ for some $i \geq 1$, $\nu_{t-} = \nu_{t_{i-1}}$.

We now introduce the probabilistic objects we need.

DEFINITION 2.4. Let $(\Omega, \mathcal{F}, P)$ be a (sufficiently large) probability space. On this space, we consider the following four independent random elements:

(i) an $\mathcal{M}$-valued random variable $\nu_0$ (the initial distribution);

(ii) a Poisson point measure $N(ds, di, dz, d\theta)$ on $[0, \infty) \times \mathbb{N}^* \times \mathbb{R}^d \times [0, 1]$, with intensity measure $\bar{\gamma} \, ds \, (\sum_{k \geq 1} \delta_k(di))(C\tilde{D}(z) \, dz) \, d\theta$ (the seed production Poisson measure);

(iii) a Poisson point measure $M(ds, di, d\theta)$ on $[0, \infty) \times \mathbb{N}^* \times [0, 1]$, with intensity measure $\bar{\mu} \, ds \, (\sum_{k \geq 1} \delta_k(di)) \, d\theta$ (the "intrinsic" death Poisson measure);

(iv) a Poisson point measure $Q(ds, di, dj, d\theta, d\theta')$ on $[0, \infty) \times \mathbb{N}^* \times \mathbb{N}^* \times [0, 1] \times [0, 1]$, with intensity measure $\bar{U}\bar{\alpha} \, ds \, (\sum_{k \geq 1} \delta_k(di))(\sum_{k \geq 1} \delta_k(dj)) \, d\theta \, d\theta'$ (the "competition" mortality Poisson measure).

We also consider the canonical filtration $(\mathcal{F}_t)_{t \geq 0}$ generated by these processes.

We finally write the BPDL model in terms of these stochastic objects.

DEFINITION 2.5. Admit Assumption A. A $(\mathcal{F}_t)_{t \geq 0}$-adapted stochastic process $\nu = (\nu_t)_{t \geq 0}$ that belongs a.s. to $\mathcal{T}$ will be called a BPDL process if a.s., for all $t \geq 0$,

$$\nu_t = \nu_0 + \int_0^t \int_{\mathbb{N}^*} \int_{\mathbb{R}^d} \int_0^1 \mathbf{1}_{\{i \leq \langle \nu_{s-}, 1 \rangle\}} \delta_{(H^i(\nu_{s-})+z)}$$

$$\times \mathbf{1}_{\{\theta \leq (\gamma(H^i(\nu_{s-}))D(H^i(\nu_{s-}),z))/(\bar{\gamma}C\tilde{D}(z))\}}$$

$$\times N(ds, di, dz, d\theta)$$

(2.8) $$- \int_0^t \int_{\mathbb{N}^*} \int_0^1 \mathbf{1}_{\{i \leq \langle \nu_{s-}, 1 \rangle\}} \delta_{H^i(\nu_{s-})} \mathbf{1}_{\{\theta \leq (\mu(H^i(\nu_{s-})))/(\bar{\mu})\}} M(ds, di, d\theta)$$

$$- \int_0^t \int_{\mathbb{N}^*} \int_{\mathbb{N}^*} \int_0^1 \int_0^1 \mathbf{1}_{\{i \leq \langle \nu_{s-}, 1 \rangle\}} \mathbf{1}_{\{j \leq \langle \nu_{s-}, 1 \rangle\}} \delta_{H^i(\nu_{s-})}$$

$$\times \mathbf{1}_{\{\theta' \leq (U(H^i(\nu_{s-}), H^j(\nu_{s-})))/(\bar{U})\}}$$

$$\times \mathbf{1}_{\{\theta \leq (\alpha(H^i(\nu_{s-})))/(\bar{\alpha})\}} Q(ds, di, dj, d\theta, d\theta').$$

Although the formula looks complicated, the principle is very simple. The indicator functions that involve $\theta$ and $\theta'$ are related to the *rates* and appear when the parameters depend on the space variable $x$. In the case where the rates are constant (studied in [2]), all the integrals and indicator functions that involve $\theta$ may be cancelled.

Let us now show that if $\nu$ solves (2.8), then it follows the dynamics in which we are interested.

PROPOSITION 2.6. *Admit Assumption* A. *Consider a solution* $(\nu_t)_{t\geq 0}$ *to* (2.8). *Then* $(\nu_t)_{t\geq 0}$ *is a Markov process. Its infinitesimal generator* $L$ *is defined for all bounded and measurable maps* $\phi : \mathcal{M} \mapsto \mathbb{R}$, *all* $\nu \in \mathcal{M}$, *by* (2.3). *In particular, the law of* $(\nu_t)_{t\geq 0}$ *does not depend on the chosen order* (*see Notation* 2.2).

PROOF. The fact that $(\nu_t)_{t\geq 0}$ is a Markov process is classical. Let us now consider a function $\phi$ as in the statement. Recall that with our notation, $\nu_0 = \sum_{i=1}^{\langle\nu_0,1\rangle} \delta_{H^i(\nu_0)}$. Recall also that $L\phi(\nu_0) = \partial_t E[\phi(\nu_t)]_{t=0}$. A simple computation, using the fact that a.s. $\phi(\nu_t) = \phi(\nu_0) + \sum_{s\leq t}[\phi(\nu_{s-} + \{\nu_s - \nu_{s-}\}) - \phi(\nu_{s-})]$, shows that

$$\phi(\nu_t) = \phi(\nu_0) + \int_0^t \int_{\mathbb{N}^*} \int_{\mathbb{R}^d} \int_0^1 [\phi(\nu_{s-} + \delta_{(H^i(\nu_{s-})+z)}) - \phi(\nu_{s-})]$$

$$\times \mathbf{1}_{\{i\leq\langle\nu_{s-},1\rangle\}}\mathbf{1}_{\{\theta\leq(\gamma(H^i(\nu_{s-}))D(H^i(\nu_{s-}),z))/(\bar{\gamma}C\tilde{D}(z))\}}$$

$$\times N(ds,di,dz,d\theta)$$

$$+ \int_0^t \int_{\mathbb{N}^*} \int_0^1 [\phi(\nu_{s-} - \delta_{H^i(\nu_{s-})}) - \phi(\nu_{s-})]$$

$$\times \mathbf{1}_{\{i\leq\langle\nu_{s-},1\rangle\}}\mathbf{1}_{\{\theta\leq(\mu(H^i(\nu_{s-})))/(\bar{\mu})\}} M(ds,di,d\theta)$$

$$+ \int_0^t \int_{\mathbb{N}^*} \int_{\mathbb{N}^*} \int_0^1 \int_0^1 [\phi(\nu_{s-} - \delta_{H^i(\nu_{s-})}) - \phi(\nu_{s-})]\mathbf{1}_{\{i\leq\langle\nu_{s-},1\rangle\}}\mathbf{1}_{\{j\leq\langle\nu_{s-},1\rangle\}}$$

$$\times \mathbf{1}_{\{\theta'\leq(U(H^i(\nu_{s-}),H^j(\nu_{s-})))/(\bar{U})\}}\mathbf{1}_{\{\theta\leq(\alpha(H^i(\nu_{s-})))/(\bar{\alpha})\}}$$

$$\times Q(ds,di,dj,d\theta,d\theta').$$

Taking expectations, we obtain

$$E[\phi(\nu_t)] = E[\phi(\nu_0)] + \int_0^t ds\, E\bigg[\int_{\mathbb{R}^d} dz\, \bar{\gamma}C\tilde{D}(z)$$

$$\times \sum_{i=1}^{\langle\nu_s,1\rangle} \frac{\gamma(H^i(\nu_{s-}))D(H^i(\nu_{s-}),z)}{\bar{\gamma}C\tilde{D}(z)}$$

$$\times [\phi(\nu_{s-} + \delta_{(H^i(\nu_{s-})+z)}) - \phi(\nu_{s-})]\Bigg]$$

$$+ \int_0^t ds\, E\Bigg[\bar{\mu} \sum_{i=1}^{\langle \nu_s, 1\rangle} \frac{\mu(H^i(\nu_{s-}))}{\bar{\mu}}[\phi(\nu_{s-} - \delta_{H^i(\nu_{s-})}) - \phi(\nu_{s-})]\Bigg]$$

$$+ \int_0^t ds\, E\Bigg[\bar{U}\bar{\alpha} \sum_{i=1}^{\langle \nu_s, 1\rangle} \sum_{j=1}^{\langle \nu_s, 1\rangle} \frac{U(H^i(\nu_{s-}), H^j(\nu_{s-}))}{\bar{U}} \frac{\alpha(H^i(\nu_{s-}))}{\bar{\alpha}}$$

$$\times [\phi(\nu_{s-} - \delta_{H^i(\nu_{s-})}) - \phi(\nu_{s-})]\Bigg]$$

$$= E[\phi(\nu_0)] + \int_0^t ds\, E\Bigg[\int_{\bar{\mathcal{X}}} \nu_s(dx) \int_{\mathbb{R}^d} dz\, \gamma(x) D(x, z)$$

$$\times [\phi(\nu_s + \delta_{(x+z)}) - \phi(\nu_s)]\Bigg]$$

$$+ \int_0^t ds\, E\Bigg[\int_{\bar{\mathcal{X}}} \nu_s(dx)[\phi(\nu_s - \delta_x) - \phi(\nu_s)]$$

$$\times \bigg\{\mu(x) + \alpha(x) \int_{\bar{\mathcal{X}}} \nu_s(dy) U(x, y)\bigg\}\Bigg].$$

Differentiating this expression at $t = 0$ leads to (2.3). □

2.3. *About simulation.* This pathwise definition of the BPDL process leads to the following simulation algorithm:

STEP 0. Simulate the initial state $\nu_0$ and set $T_0 = 0$.

STEP 1. Compute the total *event* rate, given by $m(0) = m_1(0) + m_2(0) + m_3(0)$, with

(2.9) $\quad m_1(0) = C\bar{\gamma}\langle\nu_0, 1\rangle, \qquad m_2(0) = \bar{\mu}\langle\nu_0, 1\rangle, \qquad m_3(0) = \bar{\alpha}\bar{U}\langle\nu_0, 1\rangle^2.$

Simulate $S_1$ exponentially distributed, with parameter $m(0)$, and set $T_1 = T_0 + S_1$. Set $\nu_t = \nu_0$ for all $t < T_1$. Choose whether to go to Step 1.1, 1.2 or 1.3 with probability $m_1(0)/m(0)$, $m_2(0)/m(0)$ and $m_3(0)/m(0)$.

*Step* 1.1. Choose $i$ uniformly in $\{1, \ldots, \langle\nu_0, 1\rangle\}$. Choose $z \in \mathbb{R}^d$ according to the law $\tilde{D}(z)\,dz$. With probability $1 - (\gamma(H^i(\nu_0))D(H^i(\nu_0), z))/(\bar{\gamma}C\tilde{D}(z))$, do nothing (i.e., set $\nu_{T_1} = \nu_0$); else, add a new plant at the location $H^i(\nu_0) + z$ (i.e., set $\nu_{T_1} = \nu_0 + \delta_{(H^i(\nu_0)+z)}$).

*Step* 1.2. Choose $i$ uniformly in $\{1, \ldots, \langle\nu_0, 1\rangle\}$. With probability $1 - (\mu(H^i(\nu_0)))/\bar{\mu}$, do nothing (i.e., set $\nu_{T_1} = \nu_0$); else, remove the $i$th plant (i.e., set $\nu_{T_1} = \nu_0 - \delta_{H^i(\nu_0)}$).

*Step* 1.3. Choose $i$ and $j$ uniformly in $\{1,\ldots,\langle\nu_0,1\rangle\}^2$. With probability $1-(U(H^i(\nu_0),H^j(\nu_0))\alpha(H^i(\nu_0)))/\bar{U}\bar{\alpha}$, do nothing (i.e., set $\nu_{T_1}=\nu_0$); else, remove the $i$th plant (i.e., set $\nu_{T_1}=\nu_0-\delta_{H^i(\nu_0)}$).

STEP 2. Compute the total *event* rate, given by $m(T_1)=m_1(T_1)+m_2(T_1)+m_3(T_1)$, with

(2.10)
$$m_1(T_1)=C\bar{\gamma}\langle\nu_{T_1},1\rangle,$$
$$m_2(T_1)=\bar{\mu}\langle\nu_{T_1},1\rangle,$$
$$m_3(T_1)=\bar{\alpha}\bar{U}\langle\nu_{T_1},1\rangle^2.$$

Simulate $S_2$ exponentially distributed, with parameter $m(T_1)$, and set $T_2=T_1+S_2$. Set $\nu_t=\nu_{T_1}$ for all $t\in[T_1,T_2[$ and so forth.

**3. Existence and first properties.** We now show existence, uniqueness and some moment estimates for the BPDL process.

THEOREM 3.1. (i) *Admit Assumption* A *and that* $E(\langle\nu_0,1\rangle)<\infty$. *Then there exists a unique BPDL process* $(\nu_t)_{t\geq 0}$ *in the sense of Definition* 2.5.
 (ii) *If furthermore, for some* $p\geq 1$, $E(\langle\nu_0,1\rangle^p)<\infty$, *then for any* $T<\infty$,

(3.1)
$$E\bigg(\sup_{t\in[0,T]}\langle\nu_t,1\rangle^p\bigg)<\infty.$$

PROOF. We first prove (ii). Consider thus a BPDL process $(\nu_t)_{t\geq 0}$. We introduce for each $n$ the stopping time $\tau_n=\inf\{t\geq 0,\langle\nu_t,1\rangle\geq n\}$. Then a simple computation using Assumption A shows that, neglecting the nonpositive death terms,

(3.2)
$$\sup_{s\in[0,t\wedge\tau_n]}\langle\nu_s,1\rangle^p$$
$$\leq\langle\nu_0,1\rangle^p+\int_0^{t\wedge\tau_n}\int_{\mathbb{N}^*}\int_{\mathbb{R}^d}\int_0^1[(\langle\nu_{s-},1\rangle+1)^p-\langle\nu_{s-},1\rangle^p]\mathbf{1}_{\{i\leq\langle\nu_{s-},1\rangle\}}$$
$$\times\mathbf{1}_{\{\theta\leq(\gamma(H^i(\nu_{s-}))D(H^i(\nu_{s-}),z))/(\bar{\gamma}C\tilde{D}(z))\}}$$
$$\times N(ds,di,dz,d\theta)$$
$$\leq\langle\nu_0,1\rangle^p+C_p\int_0^{t\wedge\tau_n}\int_{\mathbb{N}^*}\int_{\mathbb{R}^d}\int_0^1[1+\langle\nu_{s-},1\rangle^{p-1}]\mathbf{1}_{\{i\leq\langle\nu_{s-},1\rangle\}}$$
$$\times N(ds,di,dz,d\theta)$$

for some constant $C_p$. Taking expectations, we thus obtain, the value of $C_p$ changing from line to line:

$$
\begin{aligned}
(3.3) \quad & E\bigg(\sup_{s\in[0,t\wedge\tau_n]} \langle \nu_s, 1\rangle^p\bigg) \\
& \leq C_p + C_p E\bigg(\int_0^{t\wedge\tau_n} ds\,\bar\gamma C \int_{\mathbb{R}^d} dz\,\tilde D(z)[\langle \nu_{s-},1\rangle + \langle \nu_{s-},1\rangle^p]\bigg) \\
& \leq C_p + C_p E\bigg(\int_0^t ds\,[1 + \langle \nu_{s\wedge\tau_n},1\rangle^p]\bigg).
\end{aligned}
$$

The Gronwall lemma allows us to conclude that for any $T < \infty$, there exists a constant $C_{p,T}$, not dependent on $n$, such that

$$
(3.4) \qquad E\bigg(\sup_{t\in[0,T\wedge\tau_n]} \langle \nu_t, 1\rangle^p\bigg) \leq C_{p,T}.
$$

First, we deduce that $\tau_n$ tends a.s. to infinity. Indeed, if not, we can find a $T_0 < \infty$ such that $\varepsilon_{T_0} = P(\sup_n \tau_n < T_0) > 0$. This would imply that for all $n$, $E(\sup_{t\in[0,T_0\wedge\tau_n]}\langle \nu_t,1\rangle^p) \geq \varepsilon_{T_0} n^p$, which contradicts (3.4). We may let $n$ go to infinity in (3.4) thanks to the Fatou lemma. This leads to (3.1).

Point (i) is a consequence of point (ii). Indeed, we can, for example, build the solution $(\nu_t)_{t\geq 0}$ using the simulation algorithm previously described, and choosing the rates and acceptance–rejection according to the Poisson measures $N$, $M$ and $Q$. We have to check only that the sequence of (effective or fictitious) jump instants $T_n$ goes a.s. to infinity as $n$ tends to infinity, and this follows from (3.1) with $p = 1$. Uniqueness also holds, since we have no choice in the construction. $\square$

We now prove that if there is at most one plant at each location at time $t = 0$, then this also holds for all $t \geq 0$.

PROPOSITION 3.2. *Assume Assumption A and that $E(\langle \nu_0, 1\rangle) < \infty$. Assume also that a.s., $\sup_{x\in\bar{\mathcal{X}}} \nu_0(\{x\}) \leq 1$. Consider the Bolker–Pacala process $(\nu_t)_{t\geq 0}$. Then for all $t \geq 0$, a.s.,*

$$
(3.5) \qquad \int_{\bar{\mathcal{X}}} \nu_t(dx)\nu_t(\{x\}) = \langle \nu_t, 1\rangle \qquad \text{that is, } \sup_{x\in\bar{\mathcal{X}}} \nu_t(\{x\}) \leq 1.
$$

PROOF. Consider the nonnegative function $\phi$ defined on $\mathcal{M}$ by $\phi(\nu) = \int_{\bar{\mathcal{X}}} \nu(dx)\nu(\{x\}) - \langle \nu, 1\rangle$. Then note that a.s. $\phi(\nu_0) = 0$ and that for any $\nu \in \mathcal{M}$, any $x \in \operatorname{supp}\nu$, $\phi(\nu - \delta_x) - \phi(\nu) \leq 0$. Consider, for each $n \geq 1$, the stopping time $\tau_n = \inf\{t \geq 0, \langle \nu_t, 1\rangle \geq n\}$. A simple computation allows us to obtain, for all $t \geq 0$, all $n \geq 1$,

$$
\begin{aligned}
E[\phi(\nu_{t \wedge \tau_n})] & \\
(3.6) \qquad \leq 0 + E\bigg[\int_0^{t \wedge \tau_n} ds & \int_{\bar{\mathcal{X}}} \nu_s(dx) \int_{\mathbb{R}^d} D(x, dz) \gamma(x) \\
& \times \{\phi(\nu_s + \delta_{(x+z)}) - \phi(\nu_s)\}\bigg].
\end{aligned}
$$

We easily check, using that $\nu$ is atomic, that the right-hand side term identically vanishes, since $D(x, dz)$ has a density. Hence, a.s., $\phi(\nu_{t \wedge \tau_n}) = 0$. Thanks to (3.1) with $p = 1$, $\tau_n$ a.s. grows to infinity with $n$, which concludes the proof. □

We carry on with a property that concerns the absolute continuity of the expectation of $\nu_t$. For $\nu$ a random measure, we define the deterministic measure $E(\nu)$ by $\langle E(\nu), f \rangle = E(\langle \nu, f \rangle)$.

PROPOSITION 3.3. *Accept Assumtion A, that $E[\langle \nu_0, 1 \rangle] < \infty$ and that $E(\nu_0)$ admits a density $\tilde{n}_0$ with respect to the Lebesgue measure. Consider the BPDL process $(\nu_t)_{t \geq 0}$. Then for all $t \geq 0$, $E(\nu_t)$ has a density $\tilde{n}_t$; for all measurable nonnegative functions $f$ on $\bar{\mathcal{X}}$, $E[\langle \nu_t, f \rangle] = \int_{\bar{\mathcal{X}}} dx \, f(x) \tilde{n}_t(x)$.*

PROOF. Consider a Borel set $A$ of $\mathbb{R}^d$ with Lebesgue measure zero. Consider also, for each $n \geq 1$, the stopping time $\tau_n = \inf\{t \geq 0, \langle \nu_t, 1 \rangle \geq n\}$. A simple computation allows us to obtain, for all $t \geq 0$, all $n \geq 1$,

$$
\begin{aligned}
E[\langle \nu_{t \wedge \tau_n}, \mathbf{1}_A \rangle] = & \, E(\langle \nu_0, \mathbf{1}_A \rangle) \\
& + E\bigg(\int_0^{t \wedge \tau_n} ds \int_{\bar{\mathcal{X}}} \nu_s(dx) \gamma(x) \int_{\mathbb{R}^d} dz \, D(x, z) \mathbf{1}_A(x + z)\bigg) \\
(3.7) \qquad & - E\bigg(\int_0^{t \wedge \tau_n} ds \int_{\bar{\mathcal{X}}} \nu_s(dx) \mathbf{1}_A(x) \\
& \qquad \times \bigg(\mu(x) + \alpha(x) \int_{\bar{\mathcal{X}}} \nu_s(dy) U(x, y)\bigg)\bigg).
\end{aligned}
$$

By assumption, the first term on the right-hand side is zero. The second term is also zero, since for any $x \in \bar{\mathcal{X}}$, $\int_{\mathbb{R}^d} dz \, \mathbf{1}_A(x + z) D(x, z) = 0$. The third term is of course nonpositive. Hence for each $n$, $E(\langle \nu_{t \wedge \tau_n}, \mathbf{1}_A \rangle)$ is nonpositive and thus zero. Thanks to (3.1) with $p = 1$, $\tau_n$ a.s. grows to infinity with $n$, which concludes the proof. □

We finally give some martingale properties of the process $(\nu_t)_{t \geq 0}$.

PROPOSITION 3.4. *Admit Assumption A and that for some $p \geq 2$, $E[\langle \nu_0, 1\rangle^p] < \infty$. Consider the BPDL process $(\nu_t)_{t\geq 0}$ and recall that $L$ is defined by (2.3).*

(i) *For all measurable functions $\phi$ from $\mathcal{M}$ into $\mathbb{R}$ such that for some constant $C$, for all $\nu \in \mathcal{M}$, $|\phi(\nu)| + |L\phi(\nu)| \leq C(1 + \langle \nu, 1\rangle^p)$, the process*

$$\phi(\nu_t) - \phi(\nu_0) - \int_0^t ds\, L\phi(\nu_s) \tag{3.8}$$

*is a cadlag $L^1$-$(\mathcal{F}_t)_{t\geq 0}$-martingale starting from 0.*

(ii) *Point (i) applies to any measurable $\phi$ satisfying $|\phi(\nu)| \leq C(1+\langle \nu,1\rangle^{p-2})$.*

(iii) *Point (i) applies to any function $\phi(\nu) = \langle \nu, f \rangle^q$, with $0 \leq q \leq p-1$ and with $f$ bounded and measurable on $\bar{\mathcal{X}}$.*

(iv) *For any bounded and measurable function $f$ on $\bar{\mathcal{X}}$, the process*

$$\begin{aligned}
M_t^f = \langle \nu_t, f\rangle - \langle \nu_0, f\rangle - \int_0^t ds \int_{\bar{\mathcal{X}}} \nu_s(dx)\gamma(x) \int_{\mathbb{R}^d} dz\, D(x,z) f(x+z) \\
+ \int_0^t ds \int_{\bar{\mathcal{X}}} \nu_s(dx) f(x) \Big[ \mu(x) + \alpha(x) \int_{\bar{\mathcal{X}}} \nu_s(dy) U(x,y) \Big]
\end{aligned} \tag{3.9}$$

*is a cadlag $L^2$-martingale starting from 0 with (predictable) quadratic variation*

$$\langle M^f \rangle_t = \int_0^t ds \int_{\bar{\mathcal{X}}} \nu_s(dx) \Big\{ \gamma(x) \int_{\mathbb{R}^d} dz\, f^2(x+z) D(x,z) \\
+ f^2(x)\Big[ \mu(x) + \alpha(x) \int_{\bar{\mathcal{X}}} \nu_s(dy) U(x,y) \Big] \Big\}. \tag{3.10}$$

PROOF. First of all, note that point (i) is immediate thanks to Proposition 2.6 and (3.1). Points (ii) and (iii) follow from a straightforward computation using (2.3). To prove (iv), we first assume that $E[\langle \nu_0, 1\rangle^3] < \infty$. We apply (i) with $\phi(\nu) = \langle \nu, f\rangle$. This yields that $M^f$ is a martingale. To compute its bracket, we first apply (i) with $\phi(\nu) = \langle \nu, f\rangle^2$ and obtain that

$$\begin{aligned}
\langle \nu_t, f\rangle^2 &- \langle \nu_0, f\rangle^2 \\
&- \int_0^t ds \int_{\bar{\mathcal{X}}} \nu_s(dx)\gamma(x) \int_{\mathbb{R}^d} dz\, D(x,z) \\
&\quad \times [f^2(x+z) + 2f(x+z)\langle \nu_s, f\rangle] \\
&- \int_0^t ds \int_{\bar{\mathcal{X}}} \nu_s(dx) \{f^2(x) - 2f(x)\langle \nu_s, f\rangle\} \\
&\quad \times \Big[ \mu(x) + \alpha(x) \int_{\bar{\mathcal{X}}} \nu_s(dy) U(x,y) \Big]
\end{aligned} \tag{3.11}$$

is a martingale. Then we apply the Itô formula to compute $\langle \nu_t, f \rangle^2$ from (3.9). We deduce that

$$
\begin{aligned}
\langle \nu_t, f \rangle^2 &- \langle \nu_0, f \rangle^2 \\
&- \int_0^t ds \int_{\bar{\mathcal{X}}} \nu_s(dx) \gamma(x) \int_{\mathbb{R}^d} dz\, D(x,z) 2 f(x+z) \langle \nu_s, f \rangle \\
&+ \int_0^t ds \int_{\bar{\mathcal{X}}} \nu_s(dx) 2 f(x) \langle \nu_s, f \rangle \\
&\quad \times \left[ \mu(x) + \alpha(x) \int_{\bar{\mathcal{X}}} \nu_s(dy) U(x,y) \right] - \langle M^f \rangle_t
\end{aligned}
\tag{3.12}
$$

is a martingale. Comparing (3.11) and (3.12) leads to (3.10). The extension to the case where only $E[\langle \nu_0, 1 \rangle^2] < \infty$ is straightforward since, even in this case, $E[\langle M^f \rangle_t] < \infty$ thanks to (3.1) with $p = 2$. $\square$

**4. On the the BPDL moment equations.** We now wish to give a sense to the mean moment equation given in [2], formula (6). Note that in the biology literature, one may be confused by the notation between the discrete measure $\nu_t$, its expectation $E(\nu_t)$ [defined by $\langle E(\nu_t), f \rangle = E(\langle \nu_t, f \rangle)$] and a measure with density $n_t(x)$ of which the definition is not clear. Following [2] in this section we use the next assumption.

ASSUMPTION B. The spatial domain is $\bar{\mathcal{X}} = \mathbb{R}^d$. All the parameters $\alpha$, $\gamma$, $\mu$ and $D$ of the model are independent of $x$. Moreover, the (bounded) competition kernel $U(x,y)$ has the form $U(x - y)$, and both dispersal and competition kernels are symmetric probability distribution functions, that is, $D(z) = D(-z)$, $U(x - y) = U(y - x)$ and $\int_{\mathbb{R}^d} dz\, D(z) = \int_{\mathbb{R}^d} dz\, U(z) = 1$.

We moreover assume that $E(\langle \nu_0, 1 \rangle^2) < \infty$ and that there is at most one plant at each location at time $t = 0$. So (3.1) with $p = 1$ holds and we can define, for each time $t \in [0, T]$,

$$n(t) = E(\langle \nu_t, 1 \rangle). \tag{4.1}$$

Using Proposition 3.4(iv) with $f = 1$ and taking expectations in (3.9), we obtain

$$
\begin{aligned}
E(\langle \nu_t, 1 \rangle) &= E(\langle \nu_0, 1 \rangle) + \int_0^t ds\, (\gamma - \mu) E(\langle \nu_s, 1 \rangle) \\
&\quad - \alpha \int_0^t ds\, E\left( \int_{\mathbb{R}^d \times \mathbb{R}^d} \nu_s(dx) \nu_s(dy) U(x - y) \right).
\end{aligned}
\tag{4.2}
$$

Hence,

$$n(t) = n(0) + (\gamma - \mu)\int_0^t ds\, n(s) - \alpha \int_0^t ds\, E\bigg(\int_{\mathbb{R}^d} \nu_s(dx) U(0) \nu_s(\{x\})\bigg)$$
(4.3)
$$- \alpha \int_0^t ds\, E\bigg(\int_{\mathbb{R}^d \times \mathbb{R}^d} \nu_s(dx)\nu_s(dy)\mathbf{1}_{\{x \neq y\}} U(x-y)\bigg).$$

However, thanks to Proposition 3.2, we know that for all $s \geq 0$, $\int_{\mathbb{R}^d} \nu_s(dx) U(0) \times \nu_s(\{x\}) = U(0)\langle \nu_s, 1\rangle$. We thus obtain

$$n(t) = n(0) + (\gamma - \mu - \alpha U(0))\int_0^t ds\, n(s)$$
(4.4)
$$- \alpha \int_0^t ds\, E\bigg(\int_{\mathbb{R}^d \times \mathbb{R}^d} \nu_s(dx)\nu_s(dy)\mathbf{1}_{\{x \neq y\}} U(x-y)\bigg).$$

Let us now explain the *covariance term* used by Bolker and Pacala. Writing

$$\alpha E\bigg(\int_{\mathbb{R}^d \times \mathbb{R}^d} \nu_s(dx)\nu_s(dy)\mathbf{1}_{\{x \neq y\}} U(x-y)\bigg)$$
(4.5)
$$= \alpha E\bigg(\int_{\mathbb{R}^d \times \mathbb{R}^d} \nu_s(dx)(\nu_s(dy) - n(s)\,dy)\mathbf{1}_{\{x \neq y\}} U(x-y)\bigg) + \alpha n^2(s),$$

we obtain, from (4.4),

$$n(t) = n(0) + (\gamma - \mu - \alpha U(0))\int_0^t ds\, n(s) - \alpha \int_0^t ds\, n^2(s)$$
(4.6)
$$- \alpha \int_0^t ds\, E\bigg(\int_{\mathbb{R}^d \times \mathbb{R}^d} \nu_s(dx)(\nu_s(dy) - n(s)\,dy)\mathbf{1}_{\{x \neq y\}} U(x-y)\bigg).$$

Following the terminology of Bolker and Pacala, we define a covariance measure $C_t$ on $\mathbb{R}^d$ for each time $t$. Let $\tau_{-y}$ denote the translation by the vector $-y$. We set

(4.7) $$C_t(dr) = E\bigg(\int_{y \in \mathbb{R}^d} \mathbf{1}_{\{r \neq 0\}} \nu_t \circ \tau_{-y}^{-1}(dr) \otimes \nu_t(dy)\bigg) - n^2(t)\,dr.$$

In other words, the covariance measure is defined for each measurable bounded function $\phi$ with compact support in $\mathbb{R}^d$ by

$$\int_{\mathbb{R}^d} C_t(dr)\phi(r)$$
(4.8)
$$= E\bigg(\int_{\mathbb{R}^d \times \mathbb{R}^d} \nu_t(dx)\nu_t(dy)\mathbf{1}_{\{x \neq y\}}\phi(x-y)\bigg) - n^2(t)\int_{\mathbb{R}^d} dr\, \phi(r)$$
$$= E\bigg(\int_{\mathbb{R}^d \times \mathbb{R}^d} \nu_t(dx)(\nu_t(dy) - n(t)\,dy)\mathbf{1}_{\{x \neq y\}}\phi(x-y)\bigg).$$

By using this notation, we obtain the mean equation obtained by Bolker and Pacala ([2], formula (6), page 183), with a rigorous sense for the quadratic term:

$$(4.9) \quad \frac{dn(t)}{dt} = n(t)(\gamma - \mu - \alpha n(t)) - \alpha U(0)n(t) - \alpha \int_{\mathbb{R}^d} C_t(dr)U(r).$$

Let us finally remark that we are also able to derive an evolution equation for the covariance measure. In other words, we can write differential equations solved by $\int_{\mathbb{R}^d} C_t(dr)\phi(r)$ for all measurable bounded functions $\phi$ on $\mathbb{R}^d$ (we, however, do not obtain the same equation as in [2]). Of course moments of higher order are involved in such equations. So a remaining issue is to find reasonable *moment closures* as developed in [4]. These closures are, at the moment, not rigorously justified.

**5. Infinite particle approximations.** Our aim in this section is to describe the effect of two different normalizations on the BPDL process. In both cases, we make the initial number of plants grow to infinity. We first consider the case where the birth and death rates are unchanged. We show that the random measure $(\nu_t)_{t\geq 0}$ tends to a deterministic measure $(\xi_t)_{t\geq 0}$ and solution of a nonlinear integrodifferential equation.

In addition, the second normalization consists of accelerating the rates in a convenient way. Then $(\nu_t)_{t\geq 0}$ converges to a superprocess $(X_t)_{t\geq 0}$. This measure-valued process was introduced by Etheridge [6], who called it the *superprocess version of the Bolker–Pacala model*.

Let us first consider the most general situation.

NOTATION 5.1. For each $n \in \mathbb{N}^*$, we consider a set of parameters $(\mu_n, \gamma_n, \alpha_n, U_n, D_n)$ as in Notation 2.1, that satisfy for each $n$, Assumption A and consider an initial condition $\nu_0^n \in \mathcal{M}$. Then, we denote by $(\nu_t^n)_{t\geq 0}$ the BPDL process (see Definition 2.5) with the corresponding coefficients. We consider the subset $\mathcal{M}^n$ of $M_F(\bar{\mathcal{X}})$ defined by

$$(5.1) \quad \mathcal{M}^n = \left\{\frac{1}{n}\nu, \nu \in \mathcal{M}\right\}.$$

We finally consider the cadlag $\mathcal{M}^n$-valued Markov process $(X_t^n)_{t\geq 0}$ defined by $X_t^n = \frac{1}{n}\nu_t^n$.

The generator of $(X_t^n)_{t\geq 0}$ is then given, for any measurable map $\phi$ from $\mathcal{M}^n$ into $\mathbb{R}$, by

$$L^n\phi(\nu) = n\int_{\bar{\mathcal{X}}}\nu(dx)\int_{\mathbb{R}^d} dz\, \gamma_n(x)D_n(x,z)\left[\phi\left(\nu + \frac{1}{n}\delta_{x+z}\right) - \phi(\nu)\right]$$
(5.2)

$$+ n \int_{\bar{\mathcal{X}}} \nu(dx) \left\{ \mu_n(x) + n\alpha_n(x) \int_{\bar{\mathcal{X}}} \nu(dy) U_n(x,y) \right\}$$
$$\times \left[ \phi\left(\nu - \frac{1}{n}\delta_x\right) - \phi(\nu) \right].$$

Indeed, the generator $\tilde{L}^n$ of $(\nu_t^n)_{t\geq 0}$ is given by (2.3), replacing $(\mu, \gamma, \alpha, U, D)$ by $(\mu_n, \gamma_n, \alpha_n, U_n, D_n)$. Hence,

(5.3) $\quad L^n \phi(\nu) = \partial_t E_\nu[\phi(X_t^n)]_{t=0} = \partial_t E_{n\nu}[\phi(\nu_t^n/n)]_{t=0} = \tilde{L}^n \phi^n(n\nu),$

where $\phi^n(\mu) = \phi(\mu/n)$. The conclusion follows from a straightforward computation. We now restate Proposition 3.4 for the renormalized model.

LEMMA 5.2. *Let $n \geq 1$ be fixed and consider the process $(X_t^n)_{t \geq 0}$ defined in Notation 5.1. Assume that for some $p \geq 2$, $E[\langle X_0^n, 1 \rangle^p] < \infty$.*

(i) *For all measurable functions $\phi$ from $\mathcal{M}^n$ into $\mathbb{R}$ such that for some constant $C$, for all $\nu \in \mathcal{M}^n$, $|\phi(\nu)| + |L^n\phi(\nu)| \leq C(1 + \langle \nu, 1\rangle^p)$, the process*

(5.4) $$\phi(X_t^n) - \phi(X_0^n) - \int_0^t ds\, L^n \phi(X_s^n)$$

*is a cadlag $L^1$-martingale starting from $0$.*

(ii) *Point (i) applies to any measurable $\phi$ satisfying $|\phi(\nu)| \leq C(1+\langle \nu, 1\rangle^{p-2})$.*

(iii) *Point (i) applies to any function $\phi(\nu) = \langle \nu, f \rangle^q$, with $0 \leq q \leq p-1$ and with $f$ bounded and measurable on $\mathcal{M}$.*

(iv) *For any $f$ bounded and measurable on $\bar{\mathcal{X}}$, the process*

(5.5) $$\begin{aligned}M_t^{n,f} &= \langle X_t^n, f\rangle - \langle X_0^n, f\rangle \\ &\quad - \int_0^t ds \int_{\bar{\mathcal{X}}} X_s^n(dx) \int_{\mathbb{R}^d} dz\, \gamma_n(x) D_n(x,z) f(x+z) \\ &\quad + \int_0^t ds \int_{\bar{\mathcal{X}}} X_s^n(dx) \left\{ \mu_n(x) + n\alpha_n(x) \int_{\bar{\mathcal{X}}} X_s^n(dy) U_n(x,y) \right\} f(x)\end{aligned}$$

*is a cadlag $L^2$-martingale with (predictable) quadratic variation*

(5.6) $$\begin{aligned}\langle M^{n,f}\rangle_t &= \frac{1}{n}\int_0^t ds \int_{\bar{\mathcal{X}}} X_s^n(dx) \int_{\mathbb{R}^d} dz\, \gamma_n(x) D_n(x,z) f^2(x+z) \\ &\quad + \frac{1}{n}\int_0^t ds \int_{\bar{\mathcal{X}}} X_s^n(dx) \\ &\qquad \times \left\{ \mu_n(x) + n\alpha_n(x) \int_{\bar{\mathcal{X}}} X_s^n(dy) U_n(x,y) \right\} f^2(x).\end{aligned}$$

We endow $M_F(\bar{\mathcal{X}})$ with the weak topology.

5.1. *Convergence to a nonlinear integrodifferential equation.* Let us now consider the mean-field approximating case in which the initial number of particles $n$ tends to infinity, and the parameters of seed production and intrinsic death stay unchanged, whereas the mortality competition parameter tends to zero as $\frac{1}{n}$. We show that the BPDL process can be approximated by a deterministic nonlinear integrodifferential equation. This might be a better deterministic way to describe the model than the moment equations of [2]. In particular, it allows us to deal with space-dependent parameters.

ASSUMPTION C1.

1. The initial conditions $X_0^n$ converge in law and for the weak topology on $M_F(\bar{\mathcal{X}})$ to some deterministic finite measure $\xi_0 \in M_F(\bar{\mathcal{X}})$, and $\sup_n E(\langle X_0^n, 1 \rangle^3) < +\infty$.
2. There exist some continuous nonnegative functions $\alpha, \gamma$ and $\mu$ on $\bar{\mathcal{X}}$, bounded by $\bar{\alpha}, \bar{\gamma}$ and $\bar{\mu}$, such that $\gamma_n(x) = \gamma(x)$, $\mu_n(x) = \mu(x)$ and $\alpha_n(x) = \alpha(x)/n$.
3. There exists a bounded nonnegative symmetric continuous function $U$ on $\bar{\mathcal{X}} \times \bar{\mathcal{X}}$ bounded by $\bar{U}$ such that $U_n(x,y) = U(x,y)$.
4. There exists a continuous nonnegative function $D$ on $\bar{\mathcal{X}} \times \mathbb{R}^d$ that satisfies, for each $x \in \bar{\mathcal{X}}$, $\int_{z \in \mathbb{R}^d, x+z \in \bar{\mathcal{X}}} dz \, D(x, z) = 1$, $D(x, z) = 0$ as soon as $x + z \notin \bar{\mathcal{X}}$ and such that $D(x, z) \leq C \tilde{D}(z)$ for a constant $C > 0$ and a probability density $\tilde{D}$ on $\mathbb{R}^d$. We set $D_n(x, z) = D(x, z)$.

The first assertion of Assumption C1 is satisfied, for example, if $X_0^n = \frac{1}{n} \sum_{i=1}^n \delta_{Z^i}$, where the random variables $Z^i$ are independent, with law $\xi_0$. In this case, the number $n$ can be seen as the *volume* of particles at initial time, and the limit of $X_t^n = \frac{1}{n} \nu_t^n$ may give a rigorous sense to the *number density*.

THEOREM 5.3. *Admit Assumption* C1, *and consider the sequence of processes $X^n$ defined in Notation* 5.1. *Then for all $T > 0$, the sequence $(X^n)$ converges in law, in $\mathbb{D}([0,T], M_F(\bar{\mathcal{X}}))$, to a deterministic continuous function $(\xi_t)_{t \geq 0} \in C([0,T], M_F(\bar{\mathcal{X}}))$. This measure-valued function $\xi$ is the unique solution, satisfying $\sup_{t \in [0,T]} \langle \xi_t, 1 \rangle < \infty$, of the integrodifferential equation written in its weak form: for all bounded and measurable functions $f$ from $\bar{\mathcal{X}}$ into $\mathbb{R}$,*

$$
\begin{aligned}
\langle \xi_t, f \rangle = \langle \xi_0, f \rangle &+ \int_0^t ds \int_{\bar{\mathcal{X}}} \xi_s(dx) \gamma(x) \int_{\mathbb{R}^d} dz \, D(x,z) f(x+z) \\
&- \int_0^t ds \int_{\bar{\mathcal{X}}} \xi_s(dx) f(x) \bigg\{ \mu(x) + \alpha(x) \int_{\bar{\mathcal{X}}} \xi_s(dy) U(x,y) \bigg\}.
\end{aligned}
\tag{5.7}
$$

Note that the link between (2.8) and (5.7) is the same as the link between the continuous-time binary Galton–Watson process with birth rate $\gamma$ and death rate $\mu$, and the deterministic differential equation $f'(t) = (\gamma - \mu) f(t)$.

PROOF. We divide the proof into several steps. Let us fix $T > 0$.

STEP 3. Let us first show the uniqueness for equation (5.7). We consider two solutions $(\xi_t)_{t\geq 0}$ and $(\bar{\xi}_t)_{t\geq 0}$ of (5.7) that satisfy $\sup_{t\in[0,T]}\langle \xi_t + \bar{\xi}_t, 1\rangle = A_T < +\infty$. We consider the variation norm defined for $\mu_1$ and $\mu_2$ in $M_F(\bar{\mathcal{X}})$ by

$$\|\mu_1 - \mu_2\| = \sup_{f\in L^\infty(\bar{\mathcal{X}}),\|f\|_\infty\leq 1} |\langle \mu_1 - \mu_2, f\rangle|. \tag{5.8}$$

Then we consider some bounded and measurable function $f$ defined on $\bar{\mathcal{X}}$ such that $\|f\|_\infty \leq 1$ and we obtain

$$|\langle \xi_t - \bar{\xi}_t, f\rangle|$$

$$\leq \int_0^t ds \left| \int_{\bar{\mathcal{X}}} [\xi_s(dx) - \bar{\xi}_s(dx)] \right.$$

$$\left. \times \left(\gamma(x)\int_{\mathbb{R}^d} dz\, D(x,z) f(x+z) - \mu(x)f(x)\right) \right|$$

$$+ \int_0^t ds \left| \int_{\bar{\mathcal{X}}} [\xi_s(dx) - \bar{\xi}_s(dx)]\alpha(x)f(x)\int_{\bar{\mathcal{X}}} \xi_s(dy) U(x,y) \right| \tag{5.9}$$

$$+ \int_0^t ds \left| \int_{\bar{\mathcal{X}}} [\xi_s(dy) - \bar{\xi}_s(dy)] \int_{\bar{\mathcal{X}}} \bar{\xi}_s(dx)\alpha(x)f(x) U(x,y) \right|.$$

However, since $\|f\|_\infty \leq 1$ for all $x \in \bar{\mathcal{X}}$,

$$\left| \gamma(x)\int_{\mathbb{R}^d} dz\, D(x,z) f(x+z) - \mu(x)f(x) \right| \leq \bar\gamma + \bar\mu,$$

while

$$\left| \alpha(x)f(x)\int_{\bar{\mathcal{X}}} \xi_s(dy) U(x,y) \right| \leq \bar\alpha \bar U A_T$$

and

$$\left| \int_{\bar{\mathcal{X}}} \bar{\xi}_s(dx)\alpha(x)f(x) U(x,y) \right| \leq \bar\alpha \bar U A_T.$$

We deduce that

$$|\langle \xi_t - \bar{\xi}_t, f\rangle| \leq [\bar\gamma + \bar\mu + 2\bar\alpha\bar U A_T]\int_0^t ds\, \|\xi_s - \bar{\xi}_s\|. \tag{5.10}$$

Taking the supremum over all functions $f$ such that $\|f\|_\infty \leq 1$ and using the Gronwall lemma, we finally deduce that for all $t \leq T$, $\|\xi_t - \bar{\xi}_t\| = 0$. Uniqueness holds.

STEP 4. Let us prove some moment estimates. By Assumption C1, it is easy to check that, for all $T > 0$,

$$\sup_n E\left(\sup_{t \in [0,T]} \langle X_t^n, 1 \rangle^3 \right) < +\infty. \tag{5.11}$$

Indeed, recalling that $X_t^n = \frac{1}{n}\nu_t^n$, we can prove, following line by line the proof of Theorem 3.1(ii) with $p = 3$, that $E[\sup_{t \in [0,T]} \langle \nu_t^n, 1 \rangle^3] \leq C_T E[\langle \nu_0^n, 1 \rangle^3]$, noting that the constant $C_T$ does not depend on $n$. We easily conclude using part 1 of Assumption C1.

STEP 5. We first endow $M_F(\bar{\mathcal{X}})$ with the vague topology, the extension to the weak topology being handled in Step 6. To show the tightness of the sequence of laws $Q^n = \mathcal{L}(X^n)$ in $\mathcal{P}(\mathbb{D}([0,T], M_F(\bar{\mathcal{X}})))$, it suffices, following [15], to show that for any continuous bounded function $f$ on $\bar{\mathcal{X}}$, the sequence of laws of the processes $\langle X^n, f \rangle$ is tight in $\mathbb{D}([0,T], \mathbb{R})$. To this end, we use the Aldous criterion [1] and the Rebolledo criterion (see [7]). We have to show

$$\sup_n E\left(\sup_{t \in [0,T]} |\langle X_s^n, f \rangle| \right) < \infty, \tag{5.12}$$

and the tightness, respectively, of the laws of the martingale part and of the drift part of the semimartingales $\langle X^n, f \rangle$. Since $f$ is bounded, (5.12) is a consequence of (5.11). Let us thus consider a couple $(S, S')$ of stopping times satisfying a.s. $0 \leq S \leq S' \leq S + \delta \leq T$. Using Lemma 5.2, we get

$$\begin{aligned}
& E[|M_{S'}^{n,f} - M_S^{n,f}|] \\
& \leq E[|M_{S'}^{n,f} - M_S^{n,f}|^2]^{1/2} \leq E[\langle M^{n,f} \rangle_{S+\delta} - \langle M^{n,f} \rangle_S]^{1/2} \\
& \leq E\left[(\bar{\gamma} + \bar{\mu} + \bar{\alpha}\bar{U}) \int_S^{S+\delta} ds\, (\langle X_s^n, 1 \rangle + \langle X_s^n, 1 \rangle^2) \right]^{1/2} \leq C\sqrt{\delta},
\end{aligned} \tag{5.13}$$

where the last inequality comes from (5.11). The finite variation part of $\langle X_{S'}^n, f \rangle - \langle X_S^n, f \rangle$ is bounded by

$$\begin{aligned}
& \int_S^{S+\delta} ds\, [(\bar{\gamma} + \bar{\mu})\langle X_s^n, 1 \rangle + \bar{\alpha}\bar{U}\langle X_s^n, 1 \rangle^2] \\
& \leq \delta C\left[1 + \sup_{s \in [0,T]} \langle X_s^n, 1 \rangle^2 \right].
\end{aligned} \tag{5.14}$$

Hence, formula (5.11) allows us to conclude that the sequence $Q^n = \mathcal{L}(X^n)$ is tight.

STEP 6. Let us now denote by $Q$ the limiting law of a subsequence of $Q^n$. We still denote this subsequence by $Q^n$. Let $X = (X_t)_{t \geq 0}$ a process with law $Q$. We remark that by construction, almost surely,

$$(5.15) \qquad \sup_{t \in [0,T]} \sup_{f \in L^\infty(\bar{\mathcal{X}}), \|f\|_\infty \leq 1} |\langle X_t^n, f \rangle - \langle X_{t-}^n, f \rangle| \leq 1/n.$$

This implies that the process $X$ is a.s. strongly continuous.

STEP 7. Let us now check that a.s. the process $X$ is the unique solution of (5.7). Thanks to (5.11), it satisfies $\sup_{t \in [0,T]} \langle X_t, 1 \rangle < +\infty$ a.s. for each $T$. Standard density arguments show that it suffices to check that $X$ solves (5.7) for all $f \in C_b(\bar{\mathcal{X}})$ and all $t \geq 0$. Let thus $f \in C_b(\bar{\mathcal{X}})$ and $t \geq 0$ be fixed. For $\nu \in C([0, \infty), M_F(\bar{\mathcal{X}}))$, denote

$$(5.16) \quad \begin{aligned} \Psi_t(\nu) &= \langle \nu_t, f \rangle - \langle \nu_0, f \rangle \\ &\quad - \int_0^t ds \int_{\bar{\mathcal{X}}} \nu_s(dx) \gamma(x) \int_{\mathbb{R}^d} dz\, D(x,z) f(x+z) \\ &\quad + \int_0^t ds \int_{\bar{\mathcal{X}}} \nu_s(dx) f(x) \Big\{ \mu(x) + \alpha(x) \int_{\bar{\mathcal{X}}} \nu_s(dy) U(x,y) \Big\}. \end{aligned}$$

We have to show that

$$(5.17) \qquad E_Q[|\Psi_t(X)|] = 0.$$

However, Lemma 5.2 and Assumption C1 imply that for each $n$,

$$(5.18) \qquad M_t^{n,f} = \Psi_t(X^n).$$

A straightforward computation using Lemma 5.2, Assumption C1 and (5.11) shows that

$$(5.19) \quad E[|M_t^{n,f}|^2] = E[\langle M^{n,f} \rangle_t] \leq \frac{C_f}{n} E\bigg[ \int_0^t ds\, \{1 + \langle X_s^n, 1 \rangle^2\} \bigg] \leq \frac{C_{f,t}}{n},$$

which goes to 0 as $n$ tends to infinity. On the other hand, since $X$ is a.s. strongly continuous, since $f$ is continuous and thanks to Assumption C1, the function $\Psi_t$ is a.s. continuous at $X$. Furthermore, for any $\nu \in \mathbb{D}([0,T], M_F(\bar{\mathcal{X}}))$,

$$(5.20) \qquad |\Psi_t(\nu)| \leq C_{f,t} \sup_{s \in [0,t]} (1 + \langle \nu_s, 1 \rangle^2).$$

Hence using (5.11), we see that the sequence $(\Psi_t(X^n))_n$ is uniformly integrable and thus

$$(5.21) \qquad \lim_n E(|\Psi_t(X^n)|) = E(|\Psi_t(X)|).$$

Associating (5.18), (5.19) and (5.21), we conclude that (5.17) holds.

STEP 8. The previous steps imply that the sequence $(X^n)$ converges to $\xi$ in $\mathbb{D}([0,T], M_F(\bar{\mathcal{X}}))$, where $M_F(\bar{\mathcal{X}})$ is endowed with the vague topology. To extend the result to the case where $M_F(\bar{\mathcal{X}})$ is endowed with the weak topology, we use a criterion proved in [12]: Since the limiting process is continuous, it suffices to prove that the sequence $(\langle X^n, 1 \rangle)$ converges to $\langle \xi, 1 \rangle$ in law, in $\mathbb{D}([0,T], \bar{\mathcal{X}})$. We may of course apply Step 5 with $f \equiv 1$, which concludes the proof. $\square$

PROPOSITION 5.4. *Assume that $\xi_0$ in $M_F(\bar{\mathcal{X}})$ has a density with respect to the Lebesgue measure. Consider the associated solution $(\xi_t)_{t\geq 0}$ to (5.7). Then for every $t \geq 0$, the finite measure $\xi_t$ has a density with respect to the Lebesgue measure.*

PROOF. The proof is similar to that of Proposition 3.3. We consider a Borel subset $A$ of $\bar{\mathcal{X}}$ with measure zero. We apply (5.7) with $f = \mathbf{1}_A$. The right-hand side expression is equal to 0 since the first term is zero by hypothesis, the second one is zero since for all $x$, $\int_{\mathbb{R}^d} dz\, \mathbf{1}_{x+z \in A} D(x,z) = 0$, and the last term is nonpositive. $\square$

REMARK 5.5. (i) Equation (5.7) is the weak form of, for all $x \in \bar{\mathcal{X}}$, $t \geq 0$,

$$\partial_t \xi_t(x) = \int_{\bar{\mathcal{X}}} dy\, \xi_t(y) \gamma(y) D(y, x-y) \tag{5.22}$$
$$- \mu(x)\xi_t(x) - \alpha(x)\xi_t(x) \int_{\bar{\mathcal{X}}} dy\, \xi_t(y) U(x,y).$$

(ii) Assume now that $\bar{\mathcal{X}} = \mathbb{R}^d$, that the competition kernel is of the form $U(x,y) = U(x-y)$ and that $D(x,z) = D(z)$ does not depend on $x$. Then (5.7) is the weak form of, for all $x \in \mathbb{R}^d$, $t \geq 0$,

$$(5.23) \quad \partial_t \xi_t(x) = [\gamma \xi_t \star D](x) - \mu(x)\xi_t(x) - \alpha(x)\xi_t(x)[\xi_t \star U](x),$$

where, for example, $[\gamma \xi_t \star D](x) = \int_{\mathbb{R}^d} \xi_t(dy)\gamma(y)D(x-y)$.

5.2. *Convergence to a superprocess.* In this section we show the relationship between the original BPDL model (rigorously written in Definition 2.5) and the superprocess version of the Bolker–Pacala model introduced by Etheridge [6]. More precisely, we show that accelerating the rates of production and natural death by a factor of $n$ makes the BPDL processes converge to a continuous random measure-valued process which generalizes the one studied in [6].

ASSUMPTION C2.

1. The space $\bar{\mathcal{X}} = \mathbb{R}^d$. The initial conditions $X_0^n$ converge in law, for the weak topology on $M_F(\mathbb{R}^d)$, to a (random) measure $X_0 \in M_F(\mathbb{R}^d)$. Furthermore, $\sup_n E(\langle X_0^n, 1\rangle^3) < +\infty$.
2. There exist some continuous positive functions $\sigma(x), \alpha(x), \gamma(x)$ and $\beta(x)$ on $\mathbb{R}^d$, respectively bounded by $\bar{\sigma}, \bar{\alpha}, \bar{\gamma}$ and $\bar{\beta}$, a nonnegative symmetric continuous function $U(x,y)$ on $\mathbb{R}^d \times \mathbb{R}^d$ bounded by $\bar{U}$, such that

$$\gamma_n(x) = n\gamma(x) + \beta(x),$$
$$\mu_n(x) = n\gamma(x),$$
(5.24) $$\alpha_n(x) = \frac{\alpha(x)}{n},$$
$$U_n(x,y) = U(x,y),$$
$$D_n(x,z) = \left(\frac{n}{2\pi\sigma(x)}\right)^{d/2} \exp\left(-\frac{n|z|^2}{2\sigma(x)}\right).$$

Note that $D_n(x,z)$ is the density of a Gaussian vector with mean 0 and variance $\frac{\sigma(x)}{n}I_d$. With these coefficients and when $n$ tends to infinity, we have more and more seed production and natural death, and less and less competition. Each seed falls more and more close to its *mother*.

THEOREM 5.6. *Admit Assumption* C2 *and consider the sequence of processes* $X^n$ *defined in Notation* 5.1. *Then for all* $T > 0$, *the sequence* $(X^n)$ *converges in law, in* $\mathbb{D}([0,T], M_F(\mathbb{R}^d))$, *to the unique (in law) superprocess* $X \in C([0,T], M_F(\mathbb{R}^d))$, *defined by the conditions*

(5.25) $$\sup_{t \in [0,T]} E[\langle X_t, 1\rangle^3] < \infty$$

*and, for any* $f \in C_b^2(\mathbb{R}^d)$,

(5.26) $$\bar{M}_t^f = \langle X_t, f\rangle - \langle X_0, f\rangle - \tfrac{1}{2}\int_0^t ds \int_{\mathbb{R}^d} X_s(dx)\sigma(x)\gamma(x)\Delta f(x)$$
$$- \int_0^t ds \int_{\mathbb{R}^d} X_s(dx)f(x)\left[\beta(x) - \alpha(x)\int_{\mathbb{R}^d} X_s(dy)U(x,y)\right]$$

*is a continuous martingale with quadratic variation*

(5.27) $$\langle \bar{M}^f\rangle_t = 2\int_0^t ds \int_{\mathbb{R}^d} X_s(dx)\gamma(x)f^2(x).$$

PROOF. We break the proof into several steps.

STEP 1. Let us first prove the uniqueness of the solution of the martingale problem defined by (5.25)–(5.27); that is, the uniqueness of a probability measure $P$ on $C([0,T], M_F(\mathbb{R}^d))$ under which the canonical process $X$ satisfies (5.25)–(5.27). This result is well known for the super-Brownian process (defined by a similar martingale problem, but with $\alpha = \beta = 0$ and $\sigma = \gamma = 1$). As noted in [6], we can use the version of Dawson's Girsanov transform obtained in [5], Theorem 2.3, to deduce the uniqueness in our situation, provided the condition

$$E_P\bigg(\int_0^t ds \int_{\mathbb{R}^d} X_s(dx)\bigg[\beta(x) - \alpha(x)\int X_s(dy)U(x,y)\bigg]^2\bigg) < +\infty$$

is satisfied. This is easily obtained from (5.25) since the coefficients are bounded.

STEP 2. Next we obtain some moment estimates. First we check that for all $T < \infty$,

(5.28) $$\sup_n \sup_{t\in[0,T]} E[\langle X_t^n, 1\rangle^3] < \infty.$$

To this end, we use Lemma 5.2(i) with $\phi(\nu) = \langle \nu, 1\rangle^3$. [To be completely rigorous, first use $\phi(\nu) = \langle \nu, 1\rangle^3 \wedge A$ and then make $A$ tend to infinity.] We obtain, using Assumption C2, that for all $t \geq 0$, all $n$,

$$E\,[\langle X_t^n, 1\rangle^3]$$

(5.29)
$$= E[\langle X_0^n, 1\rangle^3] + \int_0^t ds\, E\bigg[\int_{\mathbb{R}^d} X_s^n(dx)[n^2\gamma(x) + n\beta(x)]$$
$$\times \bigg\{\bigg[\langle X_s^n, 1\rangle + \frac{1}{n}\bigg]^3 - \langle X_s^n, 1\rangle^3\bigg\}\bigg]$$
$$+ \int_0^t ds\, E\bigg[\int_{\mathbb{R}^d} X_s^n(dx)\bigg\{n^2\gamma(x) + n\alpha(x)\int_{\mathbb{R}^d} X_s^n(dy)U(x,y)\bigg\}$$
$$\times \bigg\{\bigg[\langle X_s^n, 1\rangle - \frac{1}{n}\bigg]^3 - \langle X_s^n, 1\rangle^3\bigg\}\bigg].$$

Neglecting the nonpositive competition term, we get

$$E[\langle X_t^n, 1\rangle^3]$$

$$\leq E[\langle X_0^n, 1\rangle^3]$$

(5.30)
$$+ \int_0^t ds\, E\bigg[\int_{\mathbb{R}^d} X_s^n(dx) n^2\gamma(x)$$
$$\times \bigg\{\bigg[\langle X_s^n, 1\rangle + \frac{1}{n}\bigg]^3 + \bigg[\langle X_s^n, 1\rangle - \frac{1}{n}\bigg]^3 - 2\langle X_s^n, 1\rangle^3\bigg\}\bigg]$$
$$+ \int_0^t ds\, E\bigg[\int_{\mathbb{R}^d} X_s^n(dx) n\beta(x)\bigg\{\bigg[\langle X_s^n, 1\rangle + \frac{1}{n}\bigg]^3 - \langle X_s^n, 1\rangle^3\bigg\}\bigg].$$

However, for all $x \geq 0$, all $\varepsilon \in (0,1]$, $(x+\varepsilon)^3 - x^3 \leq 6\varepsilon(1+x^2)$ and $|(x+\varepsilon)^3 + (x-\varepsilon)^3 - 2x^3| = 6\varepsilon^2 x$. We finally obtain

$$
\begin{aligned}
(5.31) \quad E[\langle X_t^n, 1\rangle^3] &\leq E[\langle X_0^n, 1\rangle^3] + 6\bar{\gamma}\int_0^t ds\, E[\langle X_s^n, 1\rangle^2] \\
&\quad + 6\bar{\beta}\int_0^t ds\, E[\langle X_s^n, 1\rangle + \langle X_s^n, 1\rangle^3].
\end{aligned}
$$

Part 1 of Assumption C2 and the Gronwall lemma allow us to conclude that (5.28) holds.

Next, we have to check that

$$(5.32) \qquad \sup_n E\left(\sup_{t\in[0,T]} \langle X_t^n, 1\rangle\right) < \infty.$$

Applying Lemma 5.2(iv) with $f \equiv 1$ and Assumption C2, we obtain

$$
\begin{aligned}
(5.33) \quad \langle X_t^n, 1\rangle &= \langle X_0^n, 1\rangle + \int_0^t ds \int_{\mathbb{R}^d} X_s^n(dx) \\
&\quad \times \left[\beta(x) - \alpha(x)\int_{\mathbb{R}^d} X_s^n(dy) U(x,y)\right] + M_t^{n,1}.
\end{aligned}
$$

Hence

$$(5.34) \qquad \sup_{s\in[0,t]} \langle X_s^n, 1\rangle \leq \langle X_0^n, 1\rangle + \bar{\beta}\int_0^t ds\, \langle X_s^n, 1\rangle + \sup_{s\in[0,t]} |M_s^{n,1}|.$$

Thanks to the Doob inequality, part 1 of Assumption C2 and the Gronwall lemma, there exists a constant $C_t$ that is not dependent on $n$ such that

$$(5.35) \qquad E\left(\sup_{s\in[0,t]} \langle X_s^n, 1\rangle\right) \leq C_t(1 + E[\langle M^{n,1}\rangle_t]^{1/2}).$$

Using (5.6) now and Assumption C2, we obtain, for some other constant $C_t$ not dependent on $n$,

$$(5.36) \quad E[\langle M^{n,1}\rangle_t] \leq (2\bar{\gamma}+\bar{\beta})\int_0^t ds\, E[\langle X_s^n, 1\rangle] + \bar{\alpha}\bar{U}\int_0^t ds\, E[\langle X_s^n, 1\rangle^2] \leq C_t$$

thanks to (5.28). This concludes the proof of (5.32).

STEP 3. We first endow $M_F(\mathbb{R}^d)$ with the vague topology. The extension to the weak topology is handled in Step 5. We prove the tightness of the sequence of laws $(\mathcal{L}(X^n))_n$ in $\mathcal{P}(\mathbb{D}([0,\infty), M_F(\mathbb{R}^d)))$ by following the same approach as in Theorem 5.3. First, we deduce from Step 2 that $\sup_n E[\sup_{s\in[0,T]}|\langle X_s^n, f\rangle|] < \infty$ for any bounded $f$. We thus have to prove that for any $f \in C_b^2(\mathbb{R}^d)$, the sequence $\langle X_t^n, f\rangle$ satisfies the Aldous–Rebolledo criterion. Let us consider a couple $(S, S')$ of stopping times satisfying a.s.

$0 \leq S \leq S' \leq S + \delta \leq T$. Using Lemma 5.2, Assumption C2 and the fact that $|\int_{\mathbb{R}^d} dz\, D_n(x,z) f(x+z) - f(x)| \leq \bar\sigma \|\Delta f\|_\infty / 2n$, we deduce the existence of a constant $C$ independent of $n$ such that the finite variation part of $\langle X^n_{S'}, f\rangle - \langle X^n_S, f\rangle$ is bounded by

$$
\begin{aligned}
&\int_S^{S+\delta} ds \int_{\mathbb{R}^d} X^n_s(dx) \bar\beta \|f\|_\infty \\
&+ \int_S^{S+\delta} ds \int_{\mathbb{R}^d} X^n_s(dx) n\gamma(x) \left| \int_{\mathbb{R}^d} dz\, D_n(x,z) f(x+z) - f(x) \right| \\
&+ \int_S^{S+\delta} ds \int_{\mathbb{R}^d} X^n_s(dx) \bar\alpha \bar U \|f\|_\infty \int_{\mathbb{R}^d} X^n_s(dy) \\
&\leq C \int_S^{S+\delta} ds\, (\langle X^n_s, 1\rangle + \langle X^n_s, 1\rangle^2).
\end{aligned}
\tag{5.37}
$$

We can also show that, for some constant $C$,

$$
E[\langle M^{n,f}\rangle_{S+\delta} - \langle M^{n,f}\rangle_S] \leq CE\left[ \int_S^{S+\delta} ds\, (\langle X^n_s, 1\rangle + \langle X^n_s, 1\rangle^2) \right]. \tag{5.38}
$$

Using the moment estimate (5.28), we finally obtain that the laws of $(M^{n,f})$ and the laws of the drift parts of $\langle X^n, f\rangle$ are tight and then, by Rebolledo's criterion, the laws of $\langle X^n, f\rangle$ are tight.

STEP 4. Let us identify the limit. Let us set $Q^n = \mathcal{L}(X^n)$, denote by $Q$ a limiting value of the tight sequence $Q^n$ and denote by $X = (X_t)_{t\geq 0}$ a process with law $Q$. Exactly as in the proof of Theorem 5.3, we can show that $X$ belongs a.s. to $C([0,T], M_F(\mathbb{R}^d))$. We have to show that $X$ satisfies conditions (5.25)–(5.27). First note that (5.25) is straightforward from (5.28). Then, we show that for any function $f$ in $C^3_b(\mathbb{R}^d)$, the process $\bar M^f_t$ defined by (5.26) is a martingale (the extension to every function in $C^2_b$ is not hard). We consider $0 \leq s_1 \leq \cdots \leq s_k < s < t$ and some continuous bounded maps $\phi_1, \ldots, \phi_k$ on $M_F(\mathbb{R}^d)$. Our aim is to prove that, if the function $\Psi$ from $\mathbb{D}([0,T], M_F(\mathbb{R}^d))$ into $\mathbb{R}$ is defined by

$$
\begin{aligned}
\Psi(\nu) &= \phi_1(\nu_{s_1}) \cdots \phi_k(\nu_{s_k}) \\
&\quad \times \Bigg\{ \langle \nu_t, f\rangle - \langle \nu_s, f\rangle - \int_s^t du\, \langle \nu_u, \gamma\sigma\, \Delta f/2\rangle \\
&\qquad - \int_s^t du \int_{\mathbb{R}^d} \nu_u(dx) f(x) \bigg[ \beta(x) - \int_{\mathbb{R}^d} \nu_u(dy) \alpha(x) U(x,y) \bigg] \Bigg\},
\end{aligned}
\tag{5.39}
$$

then

$$
E(\Psi(X)) = 0. \tag{5.40}
$$

We know from Lemma 5.2 that using Assumption C2,

$$(5.41) \quad 0 = E[\phi_1(X^n_{s_1})\cdots\phi_k(X^n_{s_k})\{M^{n,f}_t - M^{n,f}_s\}] = E[\Psi(X^n)] - A_n,$$

where $A_n$ is defined by

$$
\begin{aligned}
A_n = E\Big[ & \int_s^t du \int_{\mathbb{R}^d} X^n_u(dx) \\
& \times \Big\{\gamma(x)n\Big[\int_{\mathbb{R}^d} dz\, D_n(x,z)f(x+z) - f(x) - \frac{\sigma(x)}{2n}\Delta f(x)\Big] \\
& \quad + \beta(x)\Big[\int_{\mathbb{R}^d} dz\, D_n(x,z)f(x+z) - f(x)\Big]\Big\} \\
& \times \phi_1(X^n_{s_1})\cdots\phi_k(X^n_{s_k})\Big].
\end{aligned}
$$
(5.42)

First, an easy computation using Assumption C2, that $f$ is $C^3_b$ and (5.28) shows that

$$(5.43) \quad |A_n| \le \frac{C_f}{n}\int_s^t du\, E[\langle X^n_u, 1\rangle] \to 0$$

as $n$ grows to infinity. Next, it is clear from Assumption C2, the fact that $f$ is $C^3_b$ and that $Q$ only charges the space of continuous processes that the map $\Psi$ is $Q$-a.s. continuous. Furthermore,

$$(5.44) \quad |\Psi(\nu)| \le C\Big(1 + \langle \nu_s, 1\rangle + \langle \nu_t, 1\rangle + \int_s^t du\langle \nu_u, 1\rangle^2\Big)$$

and we easily deduce from (5.28) that the sequence $(|\Psi(X^n)|)_n$ is uniformly integrable. Hence,

$$(5.45) \quad \lim_n E(|\Psi(X^n)|) = E_Q(|\Psi(X)|).$$

Associating (5.41), (5.43) and (5.45) allows us to conclude that (5.40) holds and thus $\bar{M}^f$ is a martingale.

We finally have to show that the bracket of $\bar{M}^f$ is given by (5.27). To this end, we first check that

$$
\begin{aligned}
\bar{N}^f_t = & \langle X_t, f\rangle^2 - \langle X_0, f\rangle^2 - \int_0^t ds \int_{\mathbb{R}^d} X_s(dx) 2\gamma(x) f^2(x) \\
& - \int_0^t ds\, 2\langle X_s, f\rangle \int_{\mathbb{R}^d} X_s(dx) f(x) \\
& \times \Big[\beta(x) - \alpha(x)\int_{\mathbb{R}^d} X_s(dy) U(x,y)\Big] \\
& - \int_0^t ds\, 2\langle X_s, f\rangle \int_{\mathbb{R}^d} X_s(dx) \tfrac{1}{2}\sigma(x)\gamma(x)\Delta f(x)
\end{aligned}
$$
(5.46)

is a martingale. This can be done exactly as for $\bar{M}_t^f$, using the fact that, thanks to Lemma 5.2(iii) (with $q = 2$),

$$
\begin{aligned}
N_t^{n,f} &= \langle X_t^n, f \rangle^2 - \langle X_0^n, f \rangle^2 \\
&\quad - \int_0^t ds \int_{\mathbb{R}^d} X_s^n(dx)\gamma(x)\left[\int_{\mathbb{R}^d} dz\, f^2(x+z)D_n(x,z) + f^2(x)\right] \\
&\quad - \int_0^t ds\, 2\langle X_s^n, f\rangle \int_{\mathbb{R}^d} X_s^n(dx) \\
&\quad\quad \times \left[\beta(x)\int_{\mathbb{R}^d} dz\, f(x+z)D_n(x,z) - \alpha(x)f(x)\int_{\mathbb{R}^d} X_s^n(dy)U(x,y)\right] \\
&\quad - \int_0^t ds\, 2\langle X_s^n, f\rangle \int_{\mathbb{R}^d} X_s^n(dx)\gamma(x)n \\
&\quad\quad \times \left[\int_{\mathbb{R}^d} dz\, f(x+z)D_n(x,z) - f(x)\right] \\
&\quad - \frac{1}{n}\int_0^t ds \int_{\mathbb{R}^d} X_s^n(dx)\beta(x)\int_{\mathbb{R}^d} dz\, f^2(x+z)D_n(x,z) \\
&\quad - \frac{1}{n}\int_0^t ds \int_{\mathbb{R}^d} X_s^n(dx)\alpha(x)\int_{\mathbb{R}^d} X_s^n(dy)U(x,y)f^2(x)
\end{aligned}
\tag{5.47}
$$

is a martingale for each $n$. Next, using the Itô formula in the definition (5.26) of $\bar{M}_t^f$, we deduce that

$$
\begin{aligned}
&\langle X_t, f \rangle^2 - \langle X_0, f \rangle^2 - \langle \bar{M}^f \rangle_t \\
&\quad - \int_0^t ds\, 2\langle X_s, f\rangle \int_{\mathbb{R}^d} X_s(dx)f(x)\left[\beta(x) - \alpha(x)\int_{\mathbb{R}^d} X_s(dy)U(x,y)\right] \\
&\quad - \int_0^t ds\, 2\langle X_s, f\rangle \int_{\mathbb{R}^d} X_s(dx)\tfrac{1}{2}\sigma(x)\gamma(x)\Delta f(x)
\end{aligned}
\tag{5.48}
$$

is a martingale. Comparing this formula with (5.46) allows us to conclude that (5.27) holds.

STEP 5. The extension to the case where $M_F(\mathbb{R}^d)$ is endowed with the weak topology uses similar arguments as in Step 6 of the proof of Theorem 5.3. □

**6. About extinction and survival.** First of all, we recall a result in [6]. Consider the superprocess $X$ obtained in Theorem 5.6, and assume that $\sigma$, $\gamma$, $\beta$ and $\alpha$ are constant on $\mathbb{R}^d$. Suppose also that $U(x,y) = h(|x-y|)$ for some nonnegative decreasing function $h$ on $\mathbb{R}_+$ that satisfies $\int_0^\infty h(r)r^{d-1}\,dr < \infty$. Then if $\beta$ is sufficiently small and $\alpha$ is sufficiently large, $X$ does not survive: a.s., there exists a $t \geq 0$ such that for all $s \geq 0$, $X_{t+s} = 0$.

We can also find a complementary result in [6] which shows nonextinction with positive probability for another model—the *stepping-stone* version of the Bolker–Pacala process. Let us now come back to the BPDL process defined as the solution of (2.8). The techniques used in [6] are specific to continuous processes and cannot be generalized to the BPDL discontinuous process.

Before giving our results, let us point out the following obvious remark.

REMARK 6.1. Assume Assumption A and that $E[\langle \nu_0, 1\rangle] < \infty$. Consider the BPDL process $(\nu_t)_{t\geq 0}$. Assume also that there exist some constants $\gamma_0 \leq \mu_0$ such that for all $x \in \bar{\mathcal{X}}$, $\mu(x) \geq \mu_0$ and $\gamma(x) \leq \gamma_0$. Then $(\nu_t)_{t\geq 0}$ does a.s. not survive, that is, $P[\exists s > 0, \langle \nu_s, 1\rangle = 0] = 1$.

The proof of this remark is not hard. In such a case, the process $Z_t = \langle \nu_t, 1\rangle$ can be bounded from above by a standard continuous-time binary Galton–Watson process $Y_t$ with death rate $\mu_0$ and birth rate $\gamma_0$. Since $\mu_0 \geq \gamma_0$, extinction a.s. occurs.

In this section, we first prove almost sure extinction in a case where the state space $\bar{\mathcal{X}}$ is compact. Then we show nonextinction in the case of a discrete version of the BPDL process with a specific (and not quite realistic) competition kernel $U$.

6.1. *Extinction in the compact case.* We check a result which essentially says that if the state space $\bar{\mathcal{X}}$ is compact, then the population does almost surely not survive. Let us make the following assumption:

ASSUMPTION E.

(i) The maps $\alpha(x)$ and $\mu(x) + \alpha(x)U(x,x)$ are bounded below.

(ii) There exists a nondecreasing function $\varphi: \mathbb{R}_+ \mapsto \mathbb{R}_+$, satisfying $\varphi(0) = 0$, such that $\lim_{x\to\infty} \varphi(x) = \infty$, such that the map $x\varphi(x)$ is convex on $[0,\infty)$ and such that, for all $\nu \in \mathcal{M}$,

$$(6.1) \qquad \langle \nu \otimes \nu, U\rangle \geq \langle \nu, 1\rangle \varphi(\langle \nu, 1\rangle).$$

REMARK 6.2. Assumption E(ii) holds if $\bar{\mathcal{X}}$ is compact in $\mathbb{R}^d$, and if there exist $\varepsilon > 0$ and $\delta > 0$ such that $U(x,y) \geq \varepsilon \mathbf{1}_{\{|x-y|\leq\delta\}}$.

THEOREM 6.3. *Admit Assumptions* A *and* E, $\nu_0 \in \mathcal{M}$ *and* $E(\langle \nu_0, 1\rangle) < \infty$. *Consider the corresponding unique BPDL process* $(\nu_t)_{t\geq 0}$ *obtained in Theorem* 3.1. *Then there is almost surely extinction, that is,* $P(\exists t \geq 0, \langle \nu_t, 1\rangle = 0) = 1$.

PROOF OF REMARK 6.2. First of all, we cover $\bar{\mathcal{X}}$ with a family $\{C_l\}_{l \in \{1,\ldots,L\}}$ of disjoint cubes of $\mathbb{R}^d$ with side $\delta/\sqrt{d}$. Note that $L$ is clearly finite and that for each $l$ and each $x, y \in C_l$, $|x - y| \leq \delta$. Recall the following consequence of the Cauchy–Schwarz inequality, which says that for all $L \geq 1$ and all $\{\alpha_1, \ldots, \alpha_L\}$ in $\mathbb{R}$, $\sum_{l=1}^{L} \alpha_l^2 \geq \frac{1}{L}[\sum_{l=1}^{L} \alpha_l]^2$. Hence for all $n \geq 1$ and all $x_1, \ldots, x_n \in \bar{\mathcal{X}}$,

$$
\begin{aligned}
\sum_{i,j=1}^{n} U(x_i, x_j) &\geq \sum_{i,j=1}^{n} \varepsilon \mathbf{1}_{\{|x_i - x_j| \leq \delta\}} \geq \varepsilon \sum_{i,j=1}^{n} \sum_{l=1}^{L} \mathbf{1}_{C_l}(x_i) \mathbf{1}_{C_l}(x_j) \\
&= \varepsilon \sum_{l=1}^{L} \left[\sum_{i=1}^{n} \mathbf{1}_{C_l}(x_i)\right]^2 \geq \varepsilon \frac{1}{L} \left[\sum_{l=1}^{L} \sum_{i=1}^{n} \mathbf{1}_{C_l}(x_i)\right]^2 = \varepsilon \frac{1}{L} n^2.
\end{aligned}
\tag{6.2}
$$

We immediately deduce that for any $\nu \in \mathcal{M}$, since $\nu$ is atomic, $\langle \nu \otimes \nu, U \rangle \geq \varepsilon \frac{1}{L} \langle \nu, 1 \rangle^2$. Hence Assumption E(ii) holds with $\varphi(n) = \varepsilon \frac{1}{L} n$. $\square$

PROOF OF THEOREM 6.3. We break the proof into several steps.

STEP 6. We first of all prove that

$$A = \sup_{t \geq 0} E(\langle \nu_t, 1 \rangle) < +\infty. \tag{6.3}$$

To this end, we set $f(t) = E(\langle \nu_t, 1 \rangle)$ and use Proposition 3.4 with $\phi(\nu) = \langle \nu, 1 \rangle$ to obtain

$$f(t) = f(0) + \int_0^t ds\, E\left[\langle \nu_s, \gamma - \mu \rangle - \int_{\bar{\mathcal{X}}} \int_{\bar{\mathcal{X}}} \nu_s(dx) \nu_s(dy) \alpha(x) U(x,y)\right]. \tag{6.4}$$

Hence $f$ is differentiable. If we set $\delta = \|\gamma - \mu\|_\infty$ and $\alpha_0 = \inf_{x \in \bar{\mathcal{X}}} \alpha(x)$, we deduce that for any $t \geq 0$,

$$f'(t) \leq \delta f(t) - \alpha_0 E(\langle \nu_t \otimes \nu_t, U \rangle). \tag{6.5}$$

Using Assumption E and then the Jensen inequality, we obtain that

$$f'(t) \leq \delta f(t) - \alpha_0 f(t) \varphi(f(t)). \tag{6.6}$$

Let now $x_0$ be the greatest solution of $\delta x_0 = \alpha_0 x_0 \varphi(x_0)$ [recall that $\varphi(x)$ is nondecreasing and goes to infinity with $x$, and that $\varphi(0) = 0$]. Then we deduce from (6.6) that for any $t \geq 0$, $f(t) \leq f(0) \vee x_0$. This concludes the first step.

STEP 7. We now check that a.s.

$$\liminf_{t \to \infty} \langle \nu_t, 1 \rangle \in \{0, \infty\}. \tag{6.7}$$

Since $\langle \nu_t, 1 \rangle$ is $\mathbb{N}$-valued, it suffices to check that for any $M \in \mathbb{N}^*$,

$$P\left[\liminf_{t \to \infty} \langle \nu_t, 1 \rangle = M\right] = 0,$$

but this is clear: If $\liminf_{t \to \infty} \langle \nu_t, 1 \rangle = M$, then $\langle \nu_t, 1 \rangle$ reaches the state $M$ infinitely often, but reaches the state $M - 1$ only a finite number of times. This is (a.s.) impossible because each time $\langle \nu_t, 1 \rangle$ reaches the state $M$, the probability that its next state is $M - 1$ is bounded below by

$$(6.8) \qquad \frac{M\varepsilon_0}{M\bar{\gamma} + M\bar{\mu} + \bar{\alpha}\bar{U}M^2} > 0,$$

where $\varepsilon_0 = \inf_{x \in \bar{\mathcal{X}}}[\mu(x) + \alpha(x)U(x,x)] > 0$.

STEP 8. Since $\langle \nu_t, 1 \rangle$ is $\mathbb{N}$-valued and 0 is an absorbing state, we immediately deduce from (6.7) that a.s. $\lim_{t \to \infty} \langle \nu_t, 1 \rangle$ exists and

$$(6.9) \qquad \lim_{t \to \infty} \langle \nu_t, 1 \rangle \in \{0, \infty\}.$$

STEP 9. By Fatou's lemma and Step 1,

$$(6.10) \quad E\left[\lim_{t \to \infty} \langle \nu_t, 1 \rangle\right] = E\left[\liminf_{t \to \infty} \langle \nu_t, 1 \rangle\right] \leq \liminf_{t \to \infty} E[\langle \nu_t, 1 \rangle] \leq A.$$

Hence $\lim_{t \to \infty} \langle \nu_t, 1 \rangle < \infty$ a.s. and we deduce from (6.9) that $\lim_{t \to \infty} \langle \nu_t, 1 \rangle = 0$ a.s. This concludes the proof. □

6.2. *Survival in a simplified case.* Next, we show that in some cases, the BPDL process survives with positive probability. We are not able to handle a proof in a general case, because the problem seems very difficult. It actually looks much more difficult than the problem of survival for the contact process, which has been studied by many mathematicians (see [11]). The only result we are able to prove is deduced from a comparison with the contact process.

ASSUMPTION S.

(i) The state space $\bar{\mathcal{X}} = \mathbb{Z}^d$.
(ii) The competition kernel $U$ is pointwise, that is, $U(x,y) = \mathbf{1}_{\{x=y\}}$.
(iii) The dispersion measure $D(x, dz) = D(dz) = (1/2^d) \sum_{u \in \mathbb{Z}^d, |u|=1} \delta_u(dz)$.
(iv) $\gamma$, $\mu$ and $\alpha$ are positive constants that satisfy

$$(6.11) \qquad \frac{\gamma 2^{-d}}{\mu + \alpha} > 2.$$

Note that $\bar{\mathcal{X}} = \mathbb{Z}^d$ was not covered by our construction. The adaptation is, however, immediate.

PROPOSITION 6.4. *Admit Assumption* S, *assume that* $\nu_0 \in \mathcal{M}$, $\langle \nu_0, 1 \rangle \geq 1$ *a.s. and assume that* $E[\langle \nu_0, 1 \rangle] < \infty$. *Consider the corresponding BPDL process* $(\nu_t)_{t \geq 0}$. *This process survives with positive probability. That means that* $P(\inf_{t \geq 0} \langle \nu_t, 1 \rangle \geq 1) > 0$.

We do not handle a completely rigorous proof. To do so we would have to build a rigorous coupling between the contact process and the BPDL process.

PROOF OR PROPOSITION 6.4. We split the proof into two steps.

STEP 10. Let us first recall definitions and results about the contact process (see [11], Chapter VI). First, denote by $M_F^s$ the set of nonnegative finite measures $\eta$ on $\mathbb{Z}^d$ such that for all $x \in \mathbb{Z}^d$, $\eta(\{x\}) \in \{0, 1\}$. The contact process with parameters $\lambda_d > 0$ and $\lambda_m > 0$ is a Markov process $(\eta_t)_{t \geq 0}$, taking its values in $M_F^s$ and with generator $K$, defined for all bounded and measurable maps $\phi$ from $M_F(\mathbb{Z}^d)$ into $\mathbb{R}$ and all $\eta \in M_F(\mathbb{Z}^d)$ by

$$K\phi(\eta) = \lambda_d \int_{\mathbb{Z}^d} \eta(dx) \sum_{u \in \mathbb{Z}^d, |u|=1} \mathbf{1}_{\{\eta(\{x+u\})=0\}} [\phi(\eta + \delta_{x+u}) - \phi(\eta)]$$
(6.12)
$$+ \lambda_m \int_{\mathbb{Z}^d} \eta(dx) \mathbf{1}_{\{\eta(\{x\})=1\}} [\phi(\eta - \delta_x) - \phi(\eta)].$$

Consider an (possibly random) initial state $\eta_0$ in $M_F^s$ satisfying $\langle \eta_0, 1 \rangle \geq 1$ a.s. Then it is known (see [11], Chapter VI) that the contact process $(\eta_t)_{t \geq 0}$ with parameters $\lambda_d > 0$, $\lambda_m > 0$ and initial state $\eta_0$ exists, is unique (in law) and that under the condition $\lambda_d > 2\lambda_m$, survives with positive probability.

STEP 11. Consider now the BPDL process $(\nu_t)_{t \geq 0}$, which takes its values in the integer-valued measures on $\mathbb{Z}^d$. Denote $\tilde{\eta}_t = \sum_{x \in \mathbb{Z}^d} \mathbf{1}_{\{\nu_t(\{x\}) \geq 1\}} \delta_x$. Note that $\tilde{\eta}_t$ is always dominated by $\nu_t$. Then $(\tilde{\eta}_t)_{t \geq 0}$ is a process with values in $M_F^s$ and we can observe that $(\tilde{\eta}_t)_{t \geq 0}$ is a sort of contact process with time- and space-dependent, random parameters $\lambda_d(t, x, \omega) = \gamma 2^{-d} [1 \vee \nu_t(\{x\})]$ and $\lambda_m(t, x, \omega) = \mathbf{1}_{\nu_t(\{x\}) \leq 1} (\mu + \alpha)$. Under Assumption S, $\lambda_d(t, x, \omega)$ is uniformly bounded from below by $\underline{\lambda}_d = \gamma 2^{-d}$, while $\lambda_m(t, x, \omega)$ is uniformly bounded from above by $\overline{\lambda}_m = \mu + \alpha$. Hence, the process $(\tilde{\eta}_t)_{t \geq 0}$ is bounded below by a contact process with parameters $\underline{\lambda}_d$ and $\overline{\lambda}_m$. Since (6.11) ensures that $2\overline{\lambda}_m < \underline{\lambda}_d$, the conclusion follows from Step 1. $\square$

Note that the previously described method may not apply to the continuous-state BPDL process, since we really need the interaction to be strictly local. In fact, the only case we could treat by such a method is the case where the competition kernel is *completely local* and cannot propagate; for example, $\bar{\mathcal{X}} = \mathbb{R}^d$ and $U(x, y) \leq \sum_{p \in \mathbb{Z}^d} \mathbf{1}_{C_p}(x) \mathbf{1}_{C_p}(y)$, where, for $p \in \mathbb{Z}^d$, $C_p = [p_1, p_1 + 1] \times \cdots \times [p_d, p_d + 1]$.

**7. On equilibria.** An interesting question is that of the existence of nontrivial equilibria for the BPDL process. Since this question seems very complicated, we first try to give some results about the deterministic equation (5.7). Then we show that there exists a nontrivial equilibrium for the BPDL process that is related to the carrying capacity under a detailed balance condition which is unfortunately very restrictive. We finally present some simulations. We suppose Assumption B in the whole section.

7.1. *Equilibrium of the deterministic equation.* We first of all point out a trivial remark.

REMARK 7.1. Suppose Assumption B and that $\gamma < \mu$, and consider a nonnegative finite measure $\xi_0$ on $\mathbb{R}^d$. Consider the corresponding unique solution $(\xi_t)_{t \geq 0} \in C([0, \infty), M_F(\mathbb{R}^d))$ of (5.7). Then $\xi_t$ tends to 0 as $t$ grows to infinity in the sense that $\langle \xi_t, 1 \rangle \leq \langle \xi_0, 1 \rangle e^{-(\mu - \gamma)t}$.

This remark follows from a straightforward application of (5.7) with $f = 1$ and of the Gronwall lemma. We next generalize the existence of solutions to (5.7) to the case of possibly nonintegrable initial conditions.

PROPOSITION 7.2. *Admit Assumption* B. *Consider a nonnegative bounded measurable function $\xi_0$ on $\mathbb{R}^d$.*

1. *There exists a unique function $(\xi_t(x))_{t \geq 0, x \in \mathbb{R}^d}$ such that:*
   (i) *for all $t \geq 0$ and all $x \in \mathbb{R}^d$, $\xi_t(x) \geq 0$;*
   (ii) *for all $T < \infty$, $\sup_{t \in [0,T], x \in \mathbb{R}^d} \xi_t(x) < \infty$;*
   (iii) *for all $t \geq 0$ and all $x \in \mathbb{R}^d$,*

$$(7.1) \qquad \xi_t(x) = \xi_0(x) + \int_0^t ds\, [\gamma(\xi_s \star D)(x) - \mu \xi_s(x) - \alpha \xi_s(x)(\xi_s \star U)(x)],$$

   *where, for example, $(\xi_t \star D)(x) = \int_{\mathbb{R}^d} dy\, D(x - y) \xi_t(y)$.*
2. *For all $x \in \mathbb{R}^d$, the map $t \mapsto \xi_t(x)$ is of class $C^1$ on $[0, \infty)$, and for all $T < \infty$, $|\partial_t \xi_t(x)|$ is bounded on $[0, T] \times \mathbb{R}^d$.*
3. *If furthermore $\int_{\mathbb{R}^d} \xi_0(x)\, dx < \infty$, then for all $T < \infty$,*

$$\sup_{t \in [0,T]} \int_{\mathbb{R}^d} dx \xi_t(x) < \infty$$

   *and the finite measure-valued function $(\xi_t(x)\, dx)_{t \geq 0}$ is the unique solution to (5.7).*

Since this proposition is quite unsurprising, we only sketch the proof.

PROOF OF PROPOSITION 7.2. First note that point 2 is an immediate consequence of (7.1) and of the fact that $\xi$ is bounded, obtained in (i) and (ii). Point 3 is also easily deduced from point 1. To check the uniqueness part of point 1, it suffices to consider two solutions $(\xi_t(x))_{t\geq 0, x\in\mathbb{R}^d}$ and $(\tilde{\xi}_t(x))_{t\geq 0, x\in\mathbb{R}^d}$ to (i)–(iii), both bounded by some constant $A_T$ on $[0,T]\times\mathbb{R}^d$. A straightforward computation shows that, for $\phi(t) = \sup_{s\leq t, x\in\mathbb{R}^d}|\xi_s(x) - \tilde{\xi}_s(x)|$, for $t\leq T$,

$$\phi(t) \leq (\gamma + \mu + 2\alpha A_T)\int_0^t ds\, \phi(s). \tag{7.2}$$

[Recall that since $\int_{\mathbb{R}^d} U(x)\,dx = 1$, $\sup_{x\in\mathbb{R}^d}(\xi_s \star U)(x) \leq \sup_{x\in\mathbb{R}^d}\xi_s(x)$.] The Gronwall lemma allows us to conclude that $\xi \equiv \tilde{\xi}$.

The existence part follows from an *implicit* Picard iteration. Define $\xi_t^0(x) = \xi_0(x)$ and construct by induction a sequence of functions $(\xi_t^n)_{t\geq 0}$ such that for each $x \in \mathbb{R}^d$, $t \mapsto \xi_t^n(x)$ is of class $C^1$ on $\mathbb{R}^+$ and satisfies, for $n \geq 1$,

$$\begin{aligned}\xi_t^{n+1}(x) &= \xi_0(x) \\ &\quad + \int_0^t ds\,[\gamma(\xi_s^n \star D)(x) - \mu\xi_s^{n+1}(x) - \alpha\xi_s^{n+1}(x)(\xi_s^n \star U)(x)].\end{aligned} \tag{7.3}$$

We can, moreover, check at each step that $\xi^n$ is well defined, nonnegative and bounded on $[0,T]\times\mathbb{R}^d$ for each $n$ and each $T$. A straightforward computation shows that for all $t\geq 0$, $\sup_n \sup_{x\in\mathbb{R}^d}\xi_t^n(x) \leq \sup_{x\in\mathbb{R}^d}\xi_0(x)e^{\gamma t}$, and next that for any $T$, there exists a constant $B_T$ such that for all $t\leq T$,

$$\begin{aligned}&\sup_{x\in\mathbb{R}^d}|\xi_t^{n+1}(x) - \xi_t^n(x)| \\ &\quad\leq B_T \int_0^t ds\,\left[\sup_{x\in\mathbb{R}^d}|\xi_s^{n+1}(x) - \xi_s^n(x)| + \sup_{x\in\mathbb{R}^d}|\xi_s^n(x) - \xi_s^{n-1}(x)|\right].\end{aligned} \tag{7.4}$$

Thanks to the Gronwall lemma, we deduce that for all $T$, all $t\leq T$ and all $n$,

$$\sup_{x\in\mathbb{R}^d}|\xi_t^{n+1}(x) - \xi_t^n(x)| \leq B_T\exp(TB_T)\int_0^t ds\,\sup_{x\in\mathbb{R}^d}|\xi_s^n(x) - \xi_s^{n-1}(x)|. \tag{7.5}$$

The Picard lemma allows us to conclude that for all $T$,

$$\sum_{n\geq 1}\sup_{t\in[0,T], x\in\mathbb{R}^d}|\xi_t^{n+1}(x) - \xi_t^n(x)| < \infty. \tag{7.6}$$

Hence, there exists a function $(\xi_t(x))_{t\geq 0, x\in\mathbb{R}^d}$ such that for any $T$, $\sup_{t\in[0,T], x\in\mathbb{R}^d}|\xi_t(x) - \xi_t^n(x)|$ tends to 0. We easily check that this function satisfies points (i)–(iii). □

We may now define the equilibria.

DEFINITION 7.3. Admit Assumption B. For a nonnegative bounded continuous function $f$ on $\mathbb{R}^d$, define the function $Ff$ on $\mathbb{R}^d$ by

$$Ff(x) = \frac{\gamma[f \star D](x)}{\mu + \alpha[f \star U](x)}. \tag{7.7}$$

Then (7.1) can be rewritten as

$$\xi_t(x) = \xi_0(x) + \int_0^t ds\,(\mu + \alpha[\xi_s \star U](x))(F\xi_s(x) - \xi_s(x)). \tag{7.8}$$

This leads us to define the equilibria in the following sense. A continuous bounded nonnegative function $c$ on $\mathbb{R}^d$ is said to be a *reasonable equilibrium* of (7.1) if for all $x \in \mathbb{R}^d$,

$$c(x) = Fc(x). \tag{7.9}$$

This definition is slightly restrictive, but we may note that if $D$ and $U$ are continuous, then any solution to (7.9) such that

$$\limsup_{|x| \to \infty} [c \star D](x)/[c \star U](x) < \infty$$

will be continuous and bounded.

REMARK 7.4. Assume Assumption B, that $\gamma > \mu$ and that $\alpha > 0$. Then the constant function $c_0(x) \equiv (\gamma - \mu)/\alpha$ is a reasonable equilibrium of (7.1). The constant function $c(x) \equiv 0$ is also, of course, a reasonable equilibrium of (7.1).

Note that the quantity $(\gamma - \mu)/\alpha$ appears in [2] and is called the *carrying capacity*, which can be understood as a sort of *maximum number of plants per unit of volume*. We use the following estimate.

LEMMA 7.5. *Assume Assumption* B, *that* $\gamma > \mu$ *and that* $\alpha > 0$. *Define the signed function* $R$ *on* $\mathbb{R}^d$ *by* $R(x) = D(x) + \frac{\gamma - \mu}{\mu}(D(x) - U(x))$. *Then, for all bounded functions* $f$ *and all* $x \in \mathbb{R}^d$,

$$Ff(x) - Fc_0(x) = \frac{\mu}{\mu + \alpha[f \star U](x)}[(f - c_0) \star R](x). \tag{7.10}$$

This result is immediately proved by using simply the expression of $F$. We now state an assumption which ensures that $R(x)\,dx$ is a probability measure and hence that $F$ is a contraction around $c_0$ in the space of bounded functions.

ASSUMPTION C. $\gamma > \mu$ and for all $x \in \mathbb{R}^d$, $\gamma D(x) \geq (\gamma - \mu)U(x)$. This implies that $R(x)\,dx$ is a probability measure on $\mathbb{R}^d$.

Let us now describe a situation for which the constant function $c_0$ is the unique nontrivial reasonable equilibrium.

PROPOSITION 7.6. *Assume Assumptions* B *and* C, *that* $\gamma > 2^d \mu$ *and that* $\alpha > 0$. *Suppose also that* $D(x) = D(|x|)$, *where the map* $D$ *is nonincreasing on* $[0, \infty)$. (*This hypothesis is physically reasonable; see* [2].) *Then any nontrivial reasonable equilibrium* $c$ *of* (7.1) *identically equals* $c_0$.

PROOF. Let $c$ thus be a nontrivial reasonable equilibrium for (7.1).

STEP 1. Since $c$ is nontrivial, there exists $x_0$ such that $c(x_0) > 0$. Since $c$ is continuous, we deduce that $c$ is bounded below on a neighborhood of $x_0$. Then (7.9) and the fact that $D$ charges any neighborhood of 0 (since it is nonincreasing) ensure that $c$ never vanishes.

STEP 2. We now show that there exists a constant $\varepsilon_0 > 0$ such that for all $x \in \mathbb{R}^d$, $c(x) \geq \varepsilon_0$. To this end, we first consider $\varepsilon > 0$ such that $\gamma(1/2^d - \varepsilon) > \mu$ and then consider $a > 0$ such that $\int_{[0,a]^d} D(x)\,dx \geq 1/2^d - \varepsilon$, which is possible since $D$ is radial. Consider now any point $x = (x_1, \ldots, x_d) \in \mathbb{R}^d$ and the box $B = [x_1, x_1 + a] \times \cdots \times [x_d, x_d + a]$. Denote $m = \inf_{x \in B} c(x)$, which is positive since $c$ is continuous and never vanishes. Our aim is to show that $m \geq g(m)$, where the $C^1$ function $g$ is defined on $[0, \infty)$ by

$$g(u) = f\left[u\left(\frac{1}{2^d} - \varepsilon\right)\right],$$
(7.11)
$$f(u) = \frac{\gamma u}{\mu + \alpha\gamma/(\gamma - \mu)u}.$$

This concludes the proof of Step 2 since we can check that $g'(0) = (1/2^d - \varepsilon)\gamma/\mu > 1$ so that $m \geq \varepsilon_0 > 0$ where $\varepsilon_0$ is the smallest positive solution to $u = g(u)$.

We thus check that $m \geq g(m)$. Let $y \in B$. Using (7.9) and Assumption C, we deduce that $c(y) \geq f([c \star D](y))$. However, $f$ is nondecreasing, so that $c(y) \geq f(m \int_B dz\, D(y - z))$. Using the symmetry and the nonincreasing properties of $D$, we easily deduce that since $y \in B$, $\int_B dz\, D(y - z) \geq \int_{[0,a]^d} dz\, D(z) \geq 1/2^d - \varepsilon$. Thus for all $y \in B$, $c(y) \geq f(m(1/2^d - \varepsilon)) = g(m)$, which ends Step 2.

STEP 3. Using (7.10), Step 2 and Assumption C, we obtain

$$\sup_{x \in \mathbb{R}^d} |c(x) - c_0| = \sup_{x \in \mathbb{R}^d} |Fc(x) - Fc_0(x)|$$
(7.12)
$$\leq \frac{\mu}{\mu + \alpha\varepsilon_0} \sup_{x \in \mathbb{R}^d} |c(x) - c_0|.$$

This implies that $\sup_{x \in \mathbb{R}^d} |c(x) - c_0| = 0$. □

Although the above uniqueness result seems quite promising, we are at the moment not able to prove that under the conditions of the previous proposition, any solution $(\xi_t)_{t\geq 0}$ to (7.1) starting from a nontrivial initial condition converges to $c_0$ in some sense. We can, however, obtain two partial results.

ASSUMPTION DBC. $\alpha > 0$, $\gamma > 0$, $\mu = 0$ and $D = U$.

This assumption is a *detailed balance condition*. Indeed, under this condition, the equilibrium $c_0(x) \equiv \gamma/\alpha$ ensures that for any couple of points $x$ and $y$, the rate of appearance of plants at $x$ due to seed production at $y$ equals the rate of disappearance of plants at $x$ because of competition of plants at $y$. In other words, $\gamma D(x-y)c_0(y) = \alpha c_0(x)c_0(y)U(x-y)$. Unfortunately, this condition is very restrictive.

PROPOSITION 7.7. *Take Assumptions* B *and* DBC. *Let $\xi_0$ be a positive bounded and measurable function on $\mathbb{R}^d$. Consider the associated unique solution $(\xi_t)_{t\geq 0}$ of* (7.1) *starting from $\xi_0$ obtained in Proposition* 7.2. *Then $\xi_t$ tends to $c_0 = \gamma/\alpha$ as $t$ grows to infinity in the sense that for all $x$ and all $t$,*

$$(7.13) \qquad [\xi_t(x) - c_0]^2 \leq [\xi_0(x) - c_0]^2 \exp\left(-2\alpha[(\xi_0 \wedge c_0) \star D](x)t\right).$$

We furthermore see in the proof below that the behavior of $\xi_t$ is quite simple: If $\xi_0(x) < c_0$, then $\xi_t(x)$ increases to $c_0$, while if $\xi_0(x) > c_0$, then $\xi_t(x)$ decreases to $c_0$.

PROOF OF PROPOSITION 7.7. Since in this case, $\partial_t \xi_t(x) = -\alpha \xi_t \star D(x) \times (\xi_t(x) - c_0)$, we easily show that for all $t \geq 0$ and all $x \in \mathbb{R}^d$,

$$(7.14) \qquad \partial_t[\xi_t(x) - c_0]^2 = -2\alpha[\xi_t(x) - c_0]^2[\xi_t \star D](x).$$

Since $\xi$ is nonnegative, we deduce that $[\xi_t(x) - c_0]^2$ is nonincreasing in $t$ for each $x$. Since furthermore $\xi_t(x)$ is continuous in $t$ for each $x$, we deduce that for any $t, x$, $\xi_t(x) \geq \xi_0(x) \wedge c_0$. Hence

$$(7.15) \qquad \partial_t[\xi_t(x) - c_0]^2 \leq -2\alpha[\xi_t(x) - c_0]^2[(\xi_0 \wedge c_0) \star D](x),$$

from which the conclusion follows. □

We now treat quite a general case of coefficients $\alpha$, $\gamma$, $\mu$, $U$ and $D$, but we consider an initial condition which is only a *small perturbation* of $c_0$.

PROPOSITION 7.8. *Admit Assumptions* B *and* C, *that $\alpha > 0$ and that $U$ is bounded below by a positive continuous function $h$ on $\mathbb{R}^d$. Consider a nonnegative bounded measurable function $\xi_0$ on $\mathbb{R}^d$ such that $\int_{\mathbb{R}^d} [\xi_0(x) - c_0]^2 \, dx < \infty$. Consider the associated unique solution $(\xi_t)_{t\geq 0}$ of* (7.1) *starting from $\xi_0$ obtained in Proposition* 7.2. *Then $\xi_t$ tends to $c_0$ as $t$ grows to infinity in the sense that there exists $a > 0$ such that for all $t$,*

$$(7.16) \qquad \int_{\mathbb{R}^d} [\xi_t(x) - c_0]^2 \, dx \leq e^{-at} \int_{\mathbb{R}^d} [\xi_0(x) - c_0]^2 \, dx.$$

PROOF. We break the proof into three steps.

STEP 1. A straightforward computation using part 2 of Proposition 7.2, (7.8) and (7.10) shows that for all $t \geq 0$ and all $x \in \mathbb{R}^d$,

$$\partial_t [\xi_t(x) - c_0]^2$$
$$= 2[\xi_t(x) - c_0] \, \partial_t \xi_t(x)$$
$$= 2[\xi_t(x) - c_0][\mu + \alpha(\xi_t \star U)(x)][F\xi_t(x) - \xi_t(x)]$$
$$= 2[\xi_t(x) - c_0][\mu + \alpha(\xi_t \star U)(x)][F\xi_t(x) - Fc_0(x)]$$
$$(7.17) \qquad + 2[\xi_t(x) - c_0][\mu + \alpha(\xi_t \star U)(x)][c_0 - \xi_t(x)]$$
$$= 2\mu[\xi_t(x) - c_0][(\xi_t - c_0) \star R](x)$$
$$\quad - 2[\xi_t(x) - c_0]^2 [\mu + \alpha(\xi_t \star U)(x)]$$
$$= -2\alpha[\xi_t(x) - c_0]^2 (\xi_t \star U)(x)$$
$$\quad - 2\mu[\xi_t(x) - c_0][(\xi_t(x) - c_0) - \{(\xi_t - c_0) \star R\}(x)].$$

Integrating this differential inequality against time, we obtain

$$[\xi_t(x) - c_0]^2$$
$$(7.18) \qquad = [\xi_0(x) - c_0]^2 - 2\int_0^t ds \, \alpha[\xi_s(x) - c_0]^2 [\xi_s \star U](x)$$
$$\quad - 2\int_0^t ds \, \mu[\xi_s(x) - c_0]\{[\xi_s(x) - c_0] - [(\xi_s - c_0) \star R](x)\} \, ds.$$

Thanks to Assumption C, $R$ is a probability measure. We furthermore know that $\xi_t$, and thus $\xi_t \star U$, is bounded on $[0, T] \times \mathbb{R}^d$ for each $T$. Thus an application of the Cauchy–Schwarz and Young inequalities yields

$$(7.19) \qquad \int_{\mathbb{R}^d} dx \, [\xi_t(x) - c_0][(\xi_t(x) - c_0) \star R(x)] \leq \int_{\mathbb{R}^d} dx \, [\xi_t(x) - c_0]^2.$$

We easily deduce that for all $T \geq 0$, $\sup_{[0,T]} \int_{\mathbb{R}^d} dx \, [\xi_t(x) - c_0]^2 < \infty$. Hence (7.18) may be integrated on $x \in \mathbb{R}^d$ and we get that for all $t \geq 0$,

$$(7.20) \qquad \partial_t \int_{\mathbb{R}^d} dx \, [\xi_t(x) - c_0]^2 \leq -2\alpha \int_{\mathbb{R}^d} dx \, [\xi_t(x) - c_0]^2 [\xi_t \star U](x).$$

STEP 2. We now wish to bound $[\xi_t \star U](x)$ from below. First, we deduce from (7.20) that $\int_{\mathbb{R}^d} dx \, [\xi_t(x) - c_0]^2$ is nonincreasing in time. Hence there exists a constant $b < \infty$ such that for all $t \geq 0$,

$$(7.21) \qquad \int_{\mathbb{R}^d} dx \, \mathbf{1}_{\{\xi_t(x) \leq c_0/2\}} \leq b.$$

However, since $U(x) \geq h(x)$, for some positive continuous function $h$ there exists a constant $a > 0$ such that

$$(7.22) \qquad \inf_{A \in \mathcal{B}(\mathbb{R}^d), \int_A dx \leq b} \int_{\mathbb{R}^d/A} dz\, U(z) \geq ba.$$

Indeed, choose any compact subset $K$ of $\mathbb{R}^d$ whose Lebesgue measure equals $2b$ and set $a = \inf_{x \in K} h(x)$. Note that for all $A \in \mathcal{B}(\mathbb{R}^d)$ such that $\int_A dx \leq b$, we also have $\int_{K/A} dx \geq b$, so that

$$(7.23) \qquad \int_{\mathbb{R}^d/A} dz\, U(z) \geq \int_{K/A} dz\, h(z) \geq ba.$$

Finally using (7.22) with $A = A_{t,x} = \{y \in \mathbb{R}^d, \xi_t(x-y) \geq c_0/2\}$, of which the Lebesgue measure is smaller than $b$ thanks to (7.21), we obtain for all $x \in \mathbb{R}^d$ and all $t \geq 0$,

$$(7.24) \qquad \begin{aligned} [\xi_t \star U](x) &= \int_{\mathbb{R}^d} dy\, \xi_t(x-y) U(y) \\ &\geq \frac{c_0}{2} \int_{A_{t,x}} dy\, U(y) \geq \frac{bac_0}{2}. \end{aligned}$$

STEP 3. Gathering (7.20) and (7.24), we finally obtain

$$(7.25) \qquad \partial_t \int_{\mathbb{R}^d} dx\, [\xi_t(x) - c_0]^2 \leq -bac_0\alpha \int_{\mathbb{R}^d} dx\, [\xi_t(x) - c_0]^2$$

from which the conclusion follows. $\square$

7.2. *Equilibrium of the BPDL process.* We now to show that it might be possible to find an equilibrium for the BPDL processes. This is a first step to study the long time behavior of the BPDL process $(\nu_t)_{t \geq 0}$ defined in Definition 2.5 conditioned on nonextinction. We unfortunately are able to treat only the case where the detailed balance condition holds. Of course, such an equilibrium will be infinite. We can, however, state the following rigorous result.

We first of all denote by $\bar{\mathcal{M}}$ the set of nonnegative (possibly infinite) integer-valued measures on $\mathbb{R}^d$. We also denote by $\mathcal{A}$ the set of functions $\phi$ from $\bar{\mathcal{M}}$ into $\mathbb{R}$ of the form $\phi(\nu) = F(\langle \nu, f \rangle)$, for some bounded measurable function $F$ on $\mathbb{R}$ and some function $f$ with compact support on $\mathbb{R}^d$.

PROPOSITION 7.9. *Admit Assumptions* B *and* DBC (*see Section* 7.1) *and that* $U(0) = 0$. *Consider a Poisson measure* $\pi$ *on* $\mathbb{R}^d$ *with intensity measure* $c_0\, dx$, *where* $c_0 = \gamma/\alpha$. *Then* $\pi$ *is a stationary BPDL process in the sense that for all* $\phi \in \mathcal{A}$, $L\phi(\pi)$ *a.s. exists, belongs to* $L^1$ *and* $E[L\phi(\pi)] = 0$, *where* $L$ *is defined in* (2.3).

Note that allowing Assumption DBC and that $U(0) = 0$ implies that there is no *natural death*. We remark also that this result is somewhat surprising, since it suggests that at equilibrium, the plant locations are independent. Let us finally mention that a similar result without Assumption DBC would be much more interesting. However, the stationary process $\pi$ does not seem to be Poisson in such a case. The proof relies on the following lemma, known as Slivnyak's formula in [13] and also can be obtained from Palm measure considerations (see [8], Chapter 10).

LEMMA 7.10. *Let $\nu$ be a Poisson measure on $\mathbb{R}^d$ with intensity $m(dx)$. Denote by $\{x_i\}_{i\geq 1}$ the points of $\nu$, that is, $\nu = \sum_{i\geq 1} \delta_{x_i}$. Then for all measurable functions $h$ from $\mathbb{R}^d \times \bar{\mathcal{M}}$ into $\mathbb{R}$ such that $\int_{\mathbb{R}^d} m(dx) E[|h(x, \nu + \delta_x)|] < \infty$,*

$$(7.26) \qquad E\left[\sum_{i\geq 1} h(x_i, \nu)\right] = \int_{\mathbb{R}^d} m(dx) E[h(x, \nu + \delta_x)].$$

PROOF OF PROPOSITION 7.9. Let $\phi$ belong to $\mathcal{A}$. The fact that $L\phi(\pi)$ a.s. exists and belongs to $L^1$ for $\phi \in \mathcal{A}$ can be easily checked using the explicit expression of $L$ and standard results about Poisson measures. We thus prove only that $E[L\phi(\pi)] = 0$. Denote by $\{x_i\}_{i\geq 1}$ the points of $\pi$, that is, $\pi = \sum_{i\geq 1} \delta_{x_i}$. Hence, we obtain, using Assumption DBC,

$$\begin{aligned} E[L\phi(\pi)] &= \gamma E\left[\sum_{i\geq 1} \int_{\mathbb{R}^d} dz\, D(z)\{\phi(\pi + \delta_{x_i+z}) - \phi(\pi)\}\right] \\ (7.27) &\qquad + \alpha E\left[\sum_{i\geq 1} \{\phi(\pi - \delta_{x_i}) - \phi(\pi)\} \sum_{j\geq 1} D(x_i - x_j)\right] \\ &=: \gamma A_1 + \alpha A_2. \end{aligned}$$

We first use Lemma 7.10 with the function $h_1(x, \nu) = \int_{\mathbb{R}^d} dz\, D(z)\{\phi(\nu + \delta_{x+z}) - \phi(\nu)\}$:

$$\begin{aligned} A_1 &= E\left[\sum_{i\geq 1} h_1(x_i, \pi)\right] \\ (7.28) &= \int_{\mathbb{R}^d} c_0\, dx\, E\left[\int_{\mathbb{R}^d} dz\, D(z)\{\phi(\pi + \delta_x + \delta_{x+z}) - \phi(\pi + \delta_x)\}\right]. \end{aligned}$$

Next, with $h_2(x, \nu) = \{\phi(\nu - \delta_x) - \phi(\nu)\} \int_{\mathbb{R}^d} \nu(dy) D(x - y)$, we obtain

$$\begin{aligned} A_2 &= E\left(\sum_{i\geq 1} h_2(x_i, \pi)\right) \\ (7.29) &= \int_{\mathbb{R}^d} dx\, c_0 E\left[\{\phi(\pi) - \phi(\pi + \delta_x)\} \int_{\mathbb{R}^d} (\pi + \delta_x)(dy) D(x - y)\right]. \end{aligned}$$

Since $D(0) = U(0) = 0$, we obtain, setting $h_3^x(y, \nu) = D(x-y)\{\phi(\nu) - \phi(\nu + \delta_x)\}$,

$$A_2 = \int_{\mathbb{R}^d} dx\, c_0 E\left(\sum_{j\geq 1} h_3^x(x_j, \pi)\right). \tag{7.30}$$

Using Lemma 7.10 again, we obtain

$$\begin{aligned}
A_2 &= \int_{\mathbb{R}^d} dx\, c_0 \int_{\mathbb{R}^d} dy\, c_0 E[D(x-y)\{\phi(\pi + \delta_y) - \phi(\pi + \delta_x + \delta_y)\}] \\
&= c_0^2 \int_{\mathbb{R}^d} dx \int_{\mathbb{R}^d} dz\, E[D(z)\{\phi(\pi + \delta_x) - \phi(\pi + \delta_{x+z} + \delta_x)\}],
\end{aligned} \tag{7.31}$$

where we have used in the last equality the substitution $(y, x) \mapsto (x, x+z)$. Since $\alpha c_0^2 = \gamma c_0$, we deduce that $\gamma A_1 = -\alpha A_2$, which ends the proof. □

### 7.3. Simulations.

The previous results suggest that the BPDL process, conditioned on nonextinction, should converge as time tends to infinity to a random measure $\nu_\infty$, quite well distributed (not far from the Lebesgue measure), with $(\gamma - \mu)/\alpha$ plants per unit of volume on average. We present simulations of this situation.

We assume that $\bar{\mathcal{X}} = \mathbb{R}$ and that $\gamma = 5$, $\mu = 1$ and $\alpha = 1$. We consider the case where $U(x, y) = \mathbf{1}_{\{|x-y| \leq 1/2\}}$ and $D(z) = \frac{1}{6}\mathbf{1}_{\{|z| \leq 3\}}$. Then we compare the BPDL process $(\nu_t)_{t \geq 0}$ with the stationary solution $c_0(dx) = [(\gamma - \mu)/\alpha]\, dx$ of (7.1).

On Figure 1, we assume that $\nu_0 = \delta_0$. The boxes represent the empirical density of the BPDL process at times $t = 3$ [Figure 1(a)] and then $t = 25$ [Figure 1(b)], obtained with one simulation, while the dotted line is the density of $c_0$ [i.e., $(\gamma - \mu)/\alpha$]. We check that after some time, the BPDL process is quite well approximated by $c_0$.

Figure 2 represents the evolution in time of $\nu_t([-5, 5])$ (full line), starting either from $\nu_0 = \delta_0$ [Figure 2(a)] or from $\nu_0 = 60\delta_0$ [Figure 2(b)], and compares it with $c_0([-5, 5]) = 10(\gamma - \mu)/\alpha$ (dotted line).

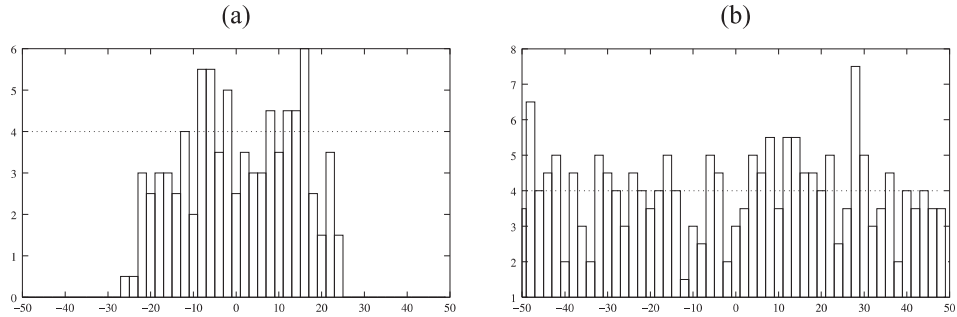

Fig. 1.  (a) $t = 3$; (b) $t = 25$.

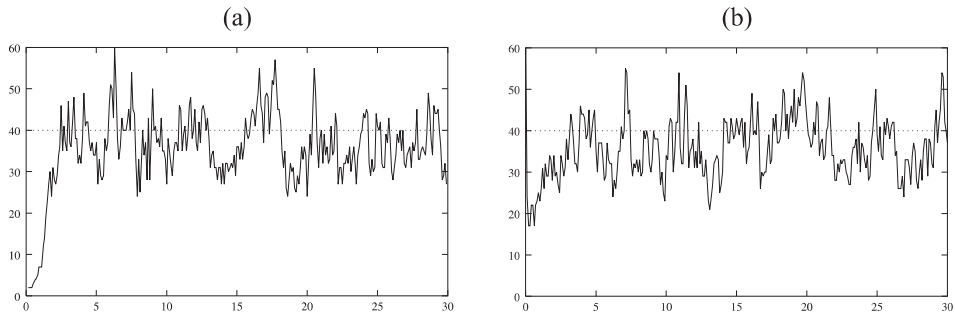

Fig. 2. (a) $\nu_0 = \delta_0$; (b) $\nu_0 = 60\delta_0$.

Finally, we measure the power of competition. To this end, we compare the evolution in time of the rate of interaction of all particles on particles located in a ball. We assume that $\nu_0 = \delta_0$. Figure 3(a) represents, in full line, the evolution in time of $\int_{\mathbb{R}} \nu_t(dx) \int_{\mathbb{R}} \nu_t(dy) \mathbf{1}_{|x| \leq 5} U(x,y)$ obtained by one simulation. The constant value (dotted line) is $\int_{\mathbb{R}} c_0(dx) \int_{\mathbb{R}} c_0(dy) \mathbf{1}_{|x| \leq 5} U(x,y) = 10 * [(\gamma - \mu)/\alpha]^2$. Figure 3(b) shows the same quantities replacing 5 by 50.

In conclusion, we can say that, on one hand, $c_0$ seems to be a good deterministic approximation of the BPDL process after a long time. On the other hand, there are clearly stochastic fluctuations around the deterministic approximation that could be interesting to study.

**Acknowledgments.** The authors thank Régis Ferrière and Bernard Roynette for numerous fruitful discussions.

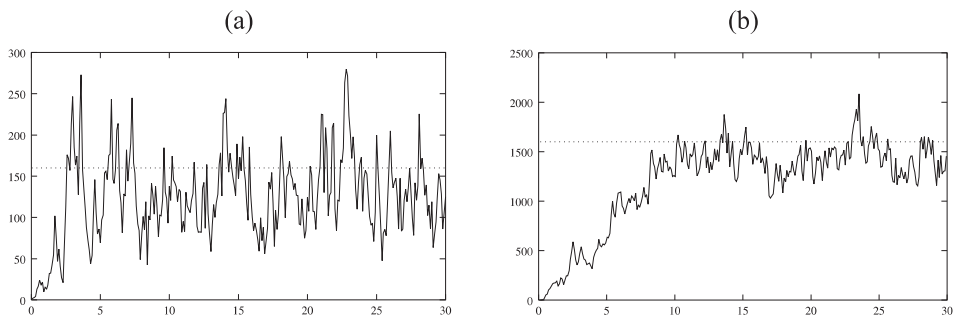

Fig. 3. *Rate of interaction endured by all particles in* (a) $[-5,5]$ *or* (b) $[-50,50]$.

FACULTÉ DES SCIENCES
INSTITUT ELIE CARTAN
BP 239
VANDOEUVRE-LÉS-NANCY CEDEX
FRANCE
E-MAIL: fournier@iecn.u-nancy.fr

MODALX
UNIVERSITÉ PARIS 10
200 AVENUE DE LA RÉPUBLIQUE
92000 NANTERRE
FRANCE
E-MAIL: sylvie.meleard@u-paris10.fr